%% file: EMCF_200316.tex
%
\newif\ifloadreferences\loadreferencestrue
%
\input preamble %
\catcode`\@=11
\def\multiline#1{\null\,\vcenter{\openup1\jot \m@th %
\ialign{\strut$\displaystyle{##}$\hfil\crcr#1\crcr}}\,}
\def\triplealign#1{\null\,\vcenter{\openup1\jot \m@th %
\ialign{\strut\hfil$\displaystyle{##}\quad$&\hfil$\displaystyle{{}##}$&$\displaystyle{{}##}$\hfil\crcr#1\crcr}}\,}
\def\tripleeqalign#1{\null\,\vcenter{\openup1\jot \m@th %
\ialign{\strut$\displaystyle{##}\quad$\hfil&$\displaystyle{{}##}$\hfil&$\displaystyle{{}##}$\hfil\crcr#1\crcr}}\,}
\catcode`\@=12
%
%
%
%
%
\myfontdefault
\global\headno=0
\global\showpagenumflag=0
\global\pageno=1
\newref{Freire}{Alikakos N. D., Freire A., The normalized mean curvature flow for a small bubble in the riemannian manifold, {\sl J. Diff. Geom.}, {\bf 64}, (2003) 247--303}
\newref{Andrews}{Andrews B., Non-collapsing in mean-convex mean curvature flow, {\sl Geometry \& Topology}, {\bf 16}, (2012), 1413--1418}
\newref{Brakke}{Brakke K. A., {\sl The motion of a surface by its mean curvature}, Princeton University Press, (1978)}
\newref{ColdMinII}{Colding T. H., Minicozzi W. P., {\sl A course in minimal surfaces}, Graduate Studies in Mathematics, {\bf 121}, American Mathematical Society, (2011)}
\newref{ColdMin}{Colding T. H., Ilmanen T., Minicozzi W. P., Rigidity of generic singularities of mean curvature flow, {\sl Publ. Math. Inst. Hautes \'Etudes Sci.}, {\bf 121}, (2015), 363--382}
\newref{GuilleminPollack}{Guillemin V., Pollack A., {\sl Differential Topology}, American Mathematical Society (Reprint), (2010)}
\newref{Hamilton}{Hamilton R., Harnack estimate for the mean curvature flow, {\sl J. Diff. Geom}, {\bf 41}, no. 1, (1995), 215--226}
\newref{Helgason}{Helgason S., {\sl Differential geometry, Lie groups and symmetric spaces}, Academic Press, (1978)}
\newref{Huisken}{Huisken G., Flow by mean curvature of convex surfaces into spheres, {\sl J. Diff. Geom.}, {\bf 20}, no. 1, (1984), 237--266}
\newref{Milnor}{Milnor J. W., {\sl Topology from the differential viewpoint}, The University Press of Virginia and Charlottesville, (1965)}
\newref{MramorPayne}{Mramor A., Payne A., Ancient and Eternal Solutions to Mean Curvature Flow from Minimal Surfaces, arXiv:1904.0843}
\newref{Salamon}{Robbin J., Salamon D., The spectral flow and the Maslov index, {\sl Bull. London Math. Soc.}, {\bf 27}, no. 1, (1995), 1--33}
\newref{SmiRos}{Rosenberg H., Smith G, Degree theory of immersed hypersurfaces, to appear in {\sl Mem. Amer. Math. Soc.}, arXiv:1010.1879}
\newref{Schwarz}{Schwarz M., {\sl Morse Homology}, Birkh\"auser Verlag, (1993)}
\newref{Smale}{Smale S., An Infinite Dimensional Version of Sard's Theorem, {\sl Amer. J. Math.}, {\bf 87}, no. 4, (1965), 861--866}
\newref{Taylor}{Taylor M. E., {\sl Partial Differential Equations I}, Applied Mathematics Sciences, {\bf 115}, Springer Verlag, (2011)}
\newref{White}{White B., The space of minimal submanifolds for varying riemannian metrics, {\sl Ind. Univ. Math. J.}, {\bf 40}, no. 1, (1991), 161--200}
\newref{WhiteII}{White B., The nature of singularities in mean curvature flow of mean-convex sets, {\sl J. Amer. Math. Soc}, {\bf 16}, no. 1, (2003), 123--138}
\makeop{a}
\def\optildea{{\text{a}}}
\makeop{cl}
\makeop{dA}
\makeop{e}

\makeop{f}
\def\optildef{{\text{f}}}
\makeop{h}
\def\opin{{\text{{\rm an}}}}
\def\optildeh{{\text{h}}}
\makeop{reg}
\makeop{s}
\makeop{C}
\makeop{CL}
\makeop{Crit}
\makeop{D}
\makeop{DI}
\makeop{E}
\makeop{GF}
\makeop{H}
\makeop{I}
\makeop{II}

\makeop{J}
\makeop{K}
\makeop{MCF}
\makeop{N}
\def\optildeN{{\tilde{\text{N}}}}
\makeop{O}
\makeop{Sign}
\makeop{Stab}
\makeop{Spec}
\makeop{T}
\makeop{V}
\makeop{W}
\makeop{Z}
\def\Pagetitle{\hfil\ifnum\pageno=1\null\else{\rm On eternal mean curvature flows...}\fi\hfil}
\def\Pagefooter{\hfil{\rm\the\pageno}\hfil}
\font\tablefont=cmr7
\def\centre{\rightskip=0pt plus 1fil \leftskip=0pt plus 1fil \spaceskip=.3333em \xspaceskip=.5em \parfillskip=0em \parindent=0em}%
\def\textmonth#1{\ifcase#1\or January\or Febuary\or March\or April\or May\or June\or July\or August\or September\or October\or November\or December\fi}
\font\abstracttitlefont=cmr10 at 12pt {\abstracttitlefont\centre On eternal mean curvature flows of tori in perturbations of the unit sphere.\par}
\medskip
{\centre Claudia V. S. Maga\~no\footnote{${}^*$}{{\tablefont Instituto de Matem\'atica, UFRJ, Av. Athos da Silveira Ramos 149, Centro de Tecnologia - Bloco C, Cidade Universit\'aria - Ilha do Fund\~ao, Caixa Postal 68530, 21941-909, Rio de Janeiro, RJ - BRAZIL\hfill}}\par}
{\centre Graham Smith${}^*$\par}
\bigskip
\noindent{\bf Abstract:~}We construct eternal mean curvature flows of tori in perturbations of the standard unit sphere $\Bbb{S}^3$. This has applications to the study of the Morse homologies of area functionals over the space of embedded tori in $\Bbb{S}^3$.
\bigskip
\newhead{Introduction}[Introduction]
\newsubhead{Eternal mean curvature flows}[EternalMeanCurvatureFlows]
Mean curvature flows, which were first studied by Brakke in \cite{Brakke}, are now well-known to geometric analysts. For the sake of clarity in what follows, we recall their formal definition. Let $\Sigma:=\Sigma^m$ be a compact, oriented, $m$-dimensional manifold. Let $M:=M^{m+1}$ be an oriented, $(m+1)$-dimensional riemannian manifold. Given an interval $I\subseteq\Bbb{R}$, which we call the {\sl time domain}, a {\sl flow} of $\Sigma$ in $M$ is defined to be a smooth function $e:\Sigma\times I\rightarrow M$ such that, for all $t$, $e_t:=e(\cdot,t)$ is an immersion. The flow $e$ is then said to be a {\sl mean curvature flow} whenever it solves the {\sl mean curvature flow equation}
$$
\bigg\langle\frac{\partial e_t}{\partial t},N_t\bigg\rangle + H_t = 0,\eqnum{\nexteqnno[MeanCurvatureFlowEquation]}
$$
where, for all $t$, $N_t$ denotes the unit normal vector field over $e_t$ compatible with the chosen orientations, and $H_t$ denotes its mean curvature with respect to this normal. This equation is invariant under the action of the {\sl reparametrisation group}, defined to be the group of all smooth functions $\alpha:\Sigma\times I\rightarrow\Sigma$ such that, for all $t$, $\alpha_t:=\alpha(\cdot,t)$ is a diffeomorphism. A mean curvature flow is said to be {\sl trivial} whenever it is constant in time, up to reparametrisation, in which case its image is trivially a minimal hypersurface in $M$. It is hard to do justice to the many interesting discoveries that have been made in the theory of mean curvature flows since Brakke's original work. For this reason, we refer the reader only to  \cite{Freire}, \cite{Andrews}, \cite{ColdMin}, \cite{Huisken}, \cite{MramorPayne} and \cite{WhiteII} for a brief selection of those results that are most relevant to the current paper.
\par
In this paper, we construct {\sl eternal mean curvature flows}. These are mean curvature flows for which the time domain coincides with the whole of $\Bbb{R}$. Since mean curvature flows tend to become singular in finite time, such flows are expected to be quite rare. With the current state of the art, it is not possible to be more precise, since little is actually known about mean curvature flows without the hypothesis of some form of convexity. What is known in the convex case consequently provides the best heuristic guide for what to expect in general. Here, in \cite{Hamilton}, Hamilton shows
\proclaim{Theorem \nextprocno, {\bf Hamilton (1995)}}
\noindent If $e$ is an eternal mean curvature flow of locally strictly convex hypersurfaces in $\Bbb{R}^{m+1}$, and if its mean curvature attains a maximum at some point of space-time, then $e$ is a translating soliton.
\endproclaim
\noindent Likewise, in \cite{WhiteII}, White conjectures
\proclaim{Conjecture \nextprocno, {\bf White (2003)}}
\noindent Every eternal, convex, non-flat mean curvature flow in $\Bbb{R}^{m+1}$ is a translating soliton.
\endproclaim
\noindent Both Hamilton's result and White's conjecture indicate that eternal mean curvature flows ought to be quite rigid. To date, few non-soliton examples of eternal mean curvature flows are known (however see \cite{MramorPayne} for a recent construction). For this reason, it is interesting to construct new examples in order to better understand their properties.
\par
In the present paper, we will only be concerned with flows of tori in perturbations of the standard unit sphere. These are of interest to us because of their relations to the work \cite{White} of White and the interesting problems suggested therein, but we see no reason why our techniques should not apply in more general settings. Let $\Bbb{T}:=\Bbb{S}^1\times\Bbb{S}^1$ be the standard $2$-dimensional torus. Let $\Bbb{S}^3$ be the standard unit sphere with constant curvature metric $g$. Let $C^\infty(\Bbb{S}^3)$ be the space of smooth functions over $\Bbb{S}^3$ furnished with the topology of smooth convergence. A subset of $C^\infty(\Bbb{S}^3)$ is said to be {\sl generic} whenever it contains a countable intersection of open, dense subsets. We prove
\proclaim{Theorem \nextprocno}
\noindent There exists a generic subset $\Cal{U}$ of $C^\infty(\Bbb{S}^3)$ with the property that, for all $u\in\Cal{U}$, there exists $\epsilon>0$ such that, for all $0<\left|t\right|<\epsilon$, there exists a non-trivial eternal mean curvature flow $e:\Bbb{T}\times\Bbb{R}\rightarrow(\Bbb{S}^3,e^{2tu}g)$.
\endproclaim
\proclabel{MainTheoremI}
\remark Theorem \procref{MainTheoremI} follows immediately from Theorems \procref{TheMorseSmaleProperty} and \procref{PerturbationTheorem}, below.
\newsubhead{The Morse-Smale property}[TheMorseSmaleProperty]
Theorem \procref{MainTheoremI} is proven via a parabolic extension of the perturbation argument developed by White in Section $3$ of \cite{White}. The main technical challenge involves the construction of functions of Morse-Smale type over the space $\opCL$ of Clifford tori in $\Bbb{S}^3$. We now describe this result, which we believe to be of independent interest.
\par
The standard {\sl Clifford torus} in $\Bbb{S}^3$ is the minimal surface
$$
T_0 := \left\{ (x,y)\ |\ \|x\|^2=\|y\|^2=1/2\right\}.\eqnum{\nexteqnno[StandardCliffordTorusIntro]}
$$
More generally, a {\sl Clifford torus} in $\Bbb{S}^3$ is any image of $T_0$ under the action of the orthogonal group $\opO(4)$. The space $\opCL$ of Clifford tori is a smooth manifold diffeomorphic to $\Bbb{R}\Bbb{P}^2\times\Bbb{R}\Bbb{P}^2$ (see Section \subheadref{TheDiffGeomOfEAndCL}). We define the linear functional $\opI:C^\infty(\Bbb{S}^3)\rightarrow C^\infty(\opCL)$ by
$$
\opI[u](T) := \int_T u\opdA_T,\eqnum{\nexteqnno[DefinitionOfFunctionalI]}
$$
where, for each $T$, $\opdA_T$ denotes its standard area form.
\par
Recall that a smooth function over any given manifold $M$ is said to be of {\sl Morse type} whenever all of its critical points are non-degenerate. Likewise, such a smooth function is said to be of {\sl Morse-Smale type} whenever, in addition to it being of Morse type, every stable manifold of its gradient flow is transverse to every unstable manifold of this flow. It is straightforward to show that the Morse-Smale property is generic amongst smooth functions over a given manifold. However, proving genericity of this property in the space of functions of the form $\opI[u]$ is far less trivial. We show
\proclaim{Theorem \nextprocno}
\noindent For generic $u\in C^\infty(\Bbb{S}^3)$, $\opI[u]$ is of Morse-Smale type.
\endproclaim
\proclabel{MainTheoremII}
\remark Theorem \procref{MainTheoremII} follows from Theorem \procref{TheMorseSmaleProperty}, below.
\newsubhead{Perturbing the Morse complex}[PerturbingTheMorseComplex]
Our results are best interpreted within the context of Morse theory. To this end, we recall the basic definitions of finite-dimensional Morse homology theory (we refer the reader to \cite{Schwarz} for a thorough introduction). Let $M:=M^m$ be a compact, $m$-dimensional riemannian manifold and let $f:M\rightarrow\Bbb{R}$ be a smooth function of Morse-Smale type. For all $k\in\Bbb{Z}$, let $\opCrit_k:=\opCrit_k(M,f)$ be the set of critical points of $f$ of Morse index equal to $k$, and define
$$
C_k := C_k(M,f) = \Bbb{Z}_2[\opCrit_k].\eqnum{\nexteqnno[DefinitionOfMorseChainGroup]}
$$
For all $k$, $C_k$ carries a natural, non-degenerate $\Bbb{Z}_2$-valued inner-product and, furthermore, $\opCrit_k$ naturally identifies with a basis of this space. For all $k$, the boundary operator $\partial_k:=\partial_k(M,f):C_k\rightarrow C_{k-1}$ is defined to be the unique linear map such that, for all $x\in\opCrit_k$ and for all $y\in\opCrit_{k-1}$, $\langle\partial_k x,y\rangle$ is equal to the number, modulo $2$, of complete gradient flows of $f$ starting at $x$ and ending at $y$, counted modulo reparametrisation. Since $f$ is of Morse-Smale type, the number of such gradient flows is always finite, so that $\partial_k$ is indeed well-defined. The complex $(C_*,\partial_*)$ is known as the {\sl Morse complex} of $(M,f)$. Since the boundary operator satisfies $\partial^2=0$, it has a well-defined homology, called the {\sl Morse homology} of $(M,f)$. This homology is independent of the Morse-Smale function $f$ chosen and is, in fact, canonically isomorphic to the singular homology of $M$.
\par
Morse homology becomes particularly interesting when applied in the infinite-di\-men\-sional setting, provided that the compactness results necessary for its development exist. In particular, let $\Cal{E}(\Bbb{T})$ be the space of smoothly embedded tori in $\Bbb{S}^3$ and, given a smooth metric $h$ over $\Bbb{S}^3$, let $\Cal{A}[h]:\Cal{E}(\Bbb{T})\rightarrow\Bbb{R}$ be its area functional. The critical points of this functional are precisely the embedded tori in $\Bbb{S}^3$ which are minimal with respect to $h$, and its complete gradient flows are precisely the eternal mean curvature flows of tori with respect to this metric. Theorem \procref{MainTheoremI} is thus reinterpreted in terms of the Morse homology of $\Cal{A}[e^{2ut}g]$ over $\Cal{E}(\Bbb{T})$ as follows.
\proclaim{Theorem \nextprocno}
\noindent For generic $u\in C^\infty(\Bbb{S}^3)$ and for sufficiently small $t$, the Morse complex of $\opI[u]$ is canonically isomorphic to a subcomplex of the Morse complex of $(\Cal{E}(\Bbb{T}),\Cal{A}[e^{2tu}g])$.
\endproclaim
\proclabel{MainTheoremIII}
\remark Theorem \procref{MainTheoremIII} follows immediately from Theorems \procref{TheMorseSmaleProperty} and \procref{PerturbationTheorem}, below.
\medskip
\remark We leave the reader to verify that the same techniques yield an analogous result when $\Cal{E}(\Bbb{T})$ is replaced with the space $\Cal{E}(\Bbb{S}^2)$ of smoothly embedded $2$-spheres in $\Bbb{S}^3$. Furthermore, since the intersection of any two equatorial spheres in $\Bbb{S}^3$ is an equatorial circle, the geometries of the intersection sets pose less of a problem, making the proof, in fact, slightly simpler (c.f. Chapter \headref{TheSpaceOfCliffordTori}).
\newhead{The space of Clifford tori}[TheSpaceOfCliffordTori]
\newsubhead{General definitions}[TheDiffGeomOfEAndCL]
We first define the geometric objects that will be used throughout the sequel. Let $\Bbb{T}:=\Bbb{S}^1\times\Bbb{S}^1$ be the standard torus. Let $\Bbb{S}^3$ be the unit sphere in $\Bbb{R}^4$. Let $T_0$ be the {\sl standard Clifford torus} in $\Bbb{S}^3$, that is
$$
T_0 := \left\{ (x,y)\ |\ \|x\|^2=\|y\|^2=1/2\right\}.\eqnum{\nexteqnno[StandardCliffordTorus]}
$$
Its {\sl standard parametrisation} is
$$
\Phi_0:\Bbb{T}\rightarrow T_0;(\theta,\phi)\mapsto\frac{1}{2}\big(\opCos(\theta),\opSin(\theta),\opCos(\phi),\opSin(\phi)\big).\eqnum{\nexteqnno[StandardParametrisation]}
$$
The {\sl Fermi parametrisation} of $\Bbb{S}^3$ about $T_0$ is
$$\eqalign{
\Phi:\Bbb{T}\times]-\pi/4,\pi/4[&\rightarrow \Bbb{S}^3;\cr
(\theta,\phi,r)&\mapsto(\opSin(\pi/4+r)\opCos(\theta),\opSin(\pi/4+r)\opSin(\theta),\cr
&\qquad\qquad\opSin(\pi/4-r)\opCos(\phi),\opSin(\pi/4-r)\opSin(\phi)),\cr}
\eqnum{\nexteqnno[FermiParametrisation]}
$$
In this parametrisation, the horizontal section $\Bbb{T}\times\left\{0\right\}$ identifies with $T_0$, the vertical lines identify with unit speed geodesics in $\Bbb{S}^3$ normal to this surface, and the riemannian metric of $\Bbb{S}^3$ is given by
$$
g = \opSin^2(\pi/4+r)d\theta^2 + \opSin^2(\pi/4-r)d\phi^2 + dr^2.\eqnum{\nexteqnno[FermiMetric]}
$$
\par
We now define the {\sl spaces} of geometric objects that will be used in the sequel. In studying infinite-dimensional manifolds, we adopt the framework of {\sl weakly smooth manifolds} developed by Rosenberg and the second author in \cite{SmiRos}. This framework allows the formal development of infinite-dimensional manifold theory with minimal technical prerequisites.
\par
Let $\hat{\Cal{E}}$ be the space of smooth embeddings $e:\Bbb{T}\rightarrow \Bbb{S}^3$ furnished with the $C^\infty$ topology. The group $\Cal{D}$ of smooth diffeomorphisms of $\Bbb{T}$ acts on the right on $\hat{\Cal{E}}$ by precomposition. The quotient space $\Cal{E}:=\hat{\Cal{E}}/\Cal{D}$, furnished with the quotient topology, is a weakly smooth manifold which we call the space of {\sl unparametrised embeddings} of $\Bbb{T}$ in $\Bbb{S}^3$. Observe that $\Cal{E}$ naturally identifies with the space $\Cal{E}'$ of those subsets of $\Bbb{S}^3$ which are smoothly embedded tori, also furnished with the topology of smooth convergence. In what follows, we will identify an element $e$ of $\hat{\Cal{E}}$, its equivalence class $[e]$ in $\Cal{E}$ and its image
$$
\opIm(e) := \left\{e(\theta,\phi)\ |\ \theta,\phi\in\Bbb{S}^1\right\}
$$
in $\Bbb{S}^3$. In particular, we consider $T_0$ as an element of $\Cal{E}$.
\par
Define $\opC:\opO(4)\rightarrow\Cal{E}$ by
$$
\opC(M) := MT_0.\eqnum{\nexteqnno[DiffeoOfOWithCL]}
$$
Its image $\opCL$ is a $4$-dimensional, strongly smooth submanifold of $\Cal{E}$ which we call the space of {\sl Clifford tori} in $\Bbb{S}^3$ (c.f. \cite{White}). To see that $\opCL$ is $4$-dimensional, observe first that $\opC$ is equivariant with respect to the left action of $\opO(4)$ on itself and the natural action of this group on $\Cal{E}$. The stabiliser of $T_0$ under the latter action is
$$
\opStab(T_0) := \bigg\{\pmatrix M\hfill&0\hfill\cr 0\hfill&N\hfill\cr\endpmatrix\ \bigg|\ M,N\in\opO(2)\bigg\},\eqnum{\nexteqnno[StabiliserOfTZero]}
$$
which naturally identifies with $\opO(2)^2$. $\opC$ therefore descends to a smooth diffeomorphism from the homogeneous space $\opO(4)/\opO(2)^2$ into $\opCL$. In particular, $\opCL$ is diffeomorphic to $\Bbb{R}\Bbb{P}^2\times\Bbb{R}\Bbb{P}^2$.
\par
Finally, the Fermi parametrisation of $\Bbb{S}^3$ defines a weakly smooth chart of $\Cal{E}$ about $T_0$ as follows. Define $\hat{\ope}:C^\infty(\Bbb{T},]-\pi/4,\pi/4[)\rightarrow\hat{\Cal{E}}$ by
$$
\hat{\ope}[f](\theta,\phi) := \Phi(\theta,\phi,f(\theta,\phi)),\eqnum{\nexteqnno[DefinitionOfHatE]}
$$
and define $\ope:C^\infty(\Bbb{T},]-\pi/4,\pi/4[)\rightarrow\Cal{E}$ by
$$
\ope := \pi\circ\hat{\ope},\eqnum{\nexteqnno[DefinitionOfE]}
$$
where $\pi:\hat{\Cal{E}}\rightarrow\Cal{E}$ here denotes the canonical projection. In other words, for all $f$, $\ope[f]$ is the image under $\Phi$ of the graph of $f$. The function $\ope$ is a weakly smooth diffeomorphism onto an open subset of $\Cal{E}$ and $\ope^{-1}(\opCL)$ is a $4$-dimensional, strongly smooth, embedded submanifold of $C^\infty(\Bbb{T},]-\pi/4,\pi/4[)$. Since the Jacobi operator of $T_0$ is given in the standard parametrisation by
$$
J := -2(\Delta+2),\eqnum{\nexteqnno[JacobiOperator]}
$$
and since all Clifford tori are minimal, it follows that the tangent space of $\ope^{-1}(\opCL)$ is given by
$$
\opT_0e^{-1}(\opCL) = \opKer(\Delta+2).\eqnum{\nexteqnno[TangentSpaceOfCL]}
$$
We will henceforth identify $\opT_{T_0}\opCL$ with $\opKer(\Delta+2)$ through $\opD\ope[0]^{-1}$.
\newsubhead{The geometry of $\opCL$}[TheSymmetricSpaceStructuresOfSAndCL]
We now study the symmetric space structure of $\opCL$. Let $\frak{o}(4)$ be the Lie algebra of $\opO(4)$ which, we recall, identifies with the space of $4$-dimensional, antisymmetric matrices. Define the positive-definite bilinear form $B$ over $\frak{o}(4)$ by
$$
B(M,N) := -\pi^2\opTr(MN),\eqnum{\nexteqnno[DefinitionOfB]}
$$
and define
$$\eqalign{
\frak{h} &:= \frak{o}(2)\oplus\frak{o}(2)\ \text{and}\cr
\frak{k} &:= \frak{h}^\perp,\cr}\eqnum{\nexteqnno[DefinitionOfHAndK]}
$$
where the orthogonal complement is taken with respect to $B$. We verify that
$$\eqalign{
[\frak{h},\frak{k}] &\subseteq \frak{k}\ \text{and}\cr
[\frak{k},\frak{k}] &\subseteq \frak{h},\cr}
$$
so that the decomposition $\frak{o}(4)=\frak{k}\oplus\frak{h}$ constitutes a polarisation of $\frak{o}(4)$. In particular $\opCL=\opO(4)/\opO(2)^2$ is a symmetric space (see \cite{Helgason}).
\par
For $A\in\opEnd(\Bbb{R}^2)$, define $\xi_A\in\frak{k}$ by
$$
\xi_A := \pmatrix 0\hfill& A\hfill\cr -A^t\hfill& 0\hfill\cr\endpmatrix.\eqnum{\nexteqnno[DefinitionOfXiA]}
$$
The map $A\mapsto\xi_A$ defines a linear isomorphism from $\opEnd(\Bbb{R}^2)$ into $\frak{k}$. Furnishing $\opEnd(\Bbb{R}^2)$ with the metric
$$
\langle A,B\rangle := 2\pi^2\opTr(AB^t)\eqnum{\nexteqnno[MetricOverEnd]}
$$
makes this map into a linear isometry. Recall now the function $\opC$ defined in Section \subheadref{TheDiffGeomOfEAndCL}. Its derivative at $\opId$ defines a linear isomorphism from $\frak{k}$ into the tangent space of $\opCL$ at $T_0$ which in turn identifies with $\opKer(\Delta+2)\subseteq C^\infty(\Bbb{T})$. For all $A\in\opEnd(\Bbb{R}^2)$, let $\psi_A$ be the image of $\xi_A$ under this map. Working through the above identifications, for all $A$, we obtain
$$
\psi_A(\theta,\phi) = (\opCos(\theta),\opSin(\theta))A\pmatrix\opCos(\phi)\hfill\cr\opSin(\phi)\hfill\cr\endpmatrix.\eqnum{\nexteqnno[FormulaForPhiA]}
$$
In particular, for all $A,B\in\opEnd(\Bbb{R}^2)$,
$$
B\big(\xi_A,\xi_B) = \langle A,B\rangle = \int_{\Bbb{T}}\psi_A\psi_B d\theta d\phi,
$$
so that $\opD\opC(\opId)$ defines a linear isometry from $\frak{k}$ into $\opT_{T_0}\opCL$, justifying the normalisations of \eqnref{DefinitionOfB} and \eqnref{MetricOverEnd}.
\par
Since $B$ is bi-invariant under the action by conjugation of $\opO(4)$ on $\frak{o}(4)$, it extends by left and right translations to a unique riemannian metric over $\opO(4)$. The tangent bundle of $\opO(4)$ then decomposes orthogonally with respect to this metric as
$$
\opT\opO(4) = \tau^L\frak{h}\oplus\tau^L\frak{k},
$$
where
$$\eqalign{
\tau^L\frak{h} &:= \left\{ (M,MA)\ |\ A\in\frak{h}\right\}\ \text{and}\cr
\tau^L\frak{k} &:= \left\{ (M,MA)\ |\ A\in\frak{k}\right\}\cr}
$$
denote respectively the left translations of $\frak{h}$ and $\frak{k}$. Since the left action of $\opO(4)$ on $\opCL$ is also isometric, it follows by the preceding paragraph that $\opC$ defines a riemannian submersion with kernel $\tau^L\frak{h}$ such that $\opC^*\opT\opCL$ is isometric to $\tau^L\frak{k}$. In particular, $\opC$ descends to a smooth isometry from $\opO(4)/\opO(2)^2$ into $\opCL$.
\newsubhead{Killing fields over $\Bbb{S}^3$}[KillingFieldsOverCL]
Let $\Phi$ be the Fermi parametrisation defined in Section \subheadref{TheDiffGeomOfEAndCL}. For all $\xi\in\frak{o}(4)$, let $X_\xi$ be the pull-back through $\Phi$ of the vector field that $\xi$ generates over $\Bbb{S}^3$ and let $F_\xi$ be its flow. Vector fields of the form $X_\xi$ are called {\sl Killing fields}. Observe that, for all $\xi$ and for all suitable $(\theta,\phi,r)$ and $t$,
$$
(\opExp(t\xi)\circ\Phi)(\theta,\phi,r) = (\Phi\circ F_{\xi,t})(\theta,\phi,r).\eqnum{\nexteqnno[FlowsAndExponentialFlows]}
$$
Since we will only require Killing fields generated by elements of $\frak{k}$, for all $A\in\opEnd(\Bbb{R}^2)$, we denote
$$\eqalign{
X_A &:= X_{\xi_A}\ \text{and}\cr
F_A &:= F_{\xi_A}.\cr}\eqnum{\nexteqnno[DefinitionOfXAndF]}
$$
\proclaim{Lemma \nextprocno}
\noindent For all $A\in\opEnd(\Bbb{R}^2)$,
$$
X_A(\theta,\phi,r) = \opCot(\pi/4+r)\frac{\partial\psi_A}{\partial\theta}(\theta,\phi)\partial_\theta - \opTan(\pi/4+r)\frac{\partial\psi_A}{\partial\phi}(\theta,\phi)\partial_\theta + \psi_A(\theta,\phi)\partial_r,\eqnum{\nexteqnno[ExplicitFormulaForX]}
$$
where $\psi_A$ is as in \eqnref{FormulaForPhiA}.
\endproclaim
\proclabel{ExplicitFormulaForX}
\proof Let $g$ be as in \eqnref{FermiMetric}. For all constant vector fields $\xi$ and $\eta$ over $\Bbb{T}\times]-\pi/4,\pi/4[$,
$$\eqalign{
(\Cal{L}_{X_A}g)(\xi,\eta)
&= \opD_{X_A}(g(\xi,\eta)) - g([X_A,\xi],\eta) - g(\xi,[X_A,\eta])\hfill\cr
&= \opD_{X_A}(g(\xi,\eta)) + g(\opD_\xi X_A,\eta) + g(\xi,\opD_{X_A}\eta).\hfill\cr}
$$
However, since $F_{A,t}$ is a flow of isometries,
$$
\Cal{L}_{X_A}g = 0.
$$
Since $g(\partial_r,\partial_r)=1$, substituting $\xi=\eta=\partial_t$ yields,
$$
g\big(\opD_{\partial_r}X_A,\partial_r\big) = -\opD_{X_A}(g(\partial_r,\partial_r)) = 0,
$$
so that the $\partial_r$ component of $X_A$ is independent of $r$. When $r=0$, this component is equal to $\psi_A$, so that
$$
X_A(\theta,\phi,r) = u(\theta,\phi,r)\partial_\theta +  v(\theta,\phi,r)\partial_\phi + \psi_A(\theta,\phi)\partial_r.
$$
Substituting $\xi=\eta=\partial_\theta$ then yields
$$
\psi_A\frac{\partial}{\partial r}\opSin^2(\pi/4+r) + 2\opSin^2(\pi/4+r)\frac{\partial u}{\partial\theta} = 0.
$$
Since $\frac{\partial^2\psi_A}{\partial\theta^2}=-\psi_A$, this becomes
$$
-\frac{\partial^2 \psi_A}{\partial\theta^2}\frac{\partial}{\partial r}\opSin^2(\pi/4+r) + 2\opSin^2(\pi/4+r)\frac{\partial u}{\partial\theta} = 0,
$$
which is solved by
$$
u(\theta,\phi,r) = \opCot(\pi/4+r)\frac{\partial\psi_A}{\partial\theta}(\theta,\phi) + \tilde{u}(\phi,r),
$$
for some function $\tilde{u}$. In the same manner, we obtain
$$
v(\theta,\phi,r) = -\opTan(\pi/4+r)\frac{\partial\psi_A}{\partial\phi}(\theta,\phi) + \tilde{v}(\theta,r),
$$
for some function $\tilde{v}$. Substituting $\xi=\partial_r$ and $\eta=\partial_\theta$ shows that $\tilde{u}$ is independent of $r$, whilst substituting $\xi=\partial_r$ and $\eta=\partial_\phi$ shows that $\tilde{v}$ is also independent of $r$. Finally, substituting $\xi=\partial_\theta$ and $\eta=\partial_\phi$ shows that
$$
\frac{\partial\tilde{v}}{\partial\theta} + \frac{\partial\tilde{u}}{\partial\phi} = 0,
$$
so that, by the principle of separation of variables,
$$
\frac{\partial\tilde{v}}{\partial\theta} = -\frac{\partial\tilde{u}}{\partial\phi} = c,
$$
for some constant $c$. Since $\Bbb{T}$ is compact, $c$ vanishes so that both $\tilde{u}$ and $\tilde{v}$ are constant. Finally, by evaluating $X_A$ at $r=0$, we show that these constants also vanish. This completes the proof.\qed
\newsubhead{The geometry of curves in $\opCL$}[TheGeometryOfCurvesInCL]
We now derive an asymptotic formula for curves in $\opCL$ passing through $T_0$ with prescribed first and second derivatives at this point. Recall from Section \subheadref{TheDiffGeomOfEAndCL} that $\opT_{T_0}\opCL$ identifies with $\opKer(\Delta+2)$. Choose $A,B\in\opEnd(\Bbb{R}^2)$ and let $\gamma:\Bbb{R}\rightarrow\opCL$ be such that
$$\eqalign{
\gamma(0) &= T_0,\hfill\cr
\dot{\gamma}(0) &= \psi_A\ \text{and}\hfill\cr
(\nabla_{\dot{\gamma}}\dot{\gamma})(0) &= \psi_B,\hfill\cr}\eqnum{\nexteqnno[DefinitionOfGammaInCL]}
$$
where $\nabla$ here denotes the Levi-Civita covariant derivative of $\opCL$. For sufficiently small $\epsilon$, let $f:]-\epsilon,\epsilon[\rightarrow C^\infty(\Bbb{T},]-\pi/4,\pi/4[)$ be such that, for all $t$,
$$
\ope[f_t] = \gamma(t).\eqnum{\nexteqnno[DefinitionOfFtInCL]}
$$
\proclaim{Lemma \nextprocno}
\noindent If $A$, $B$ and $\opf$ are as above, then
$$
f_t = t\psi_A + \frac{1}{2}t^2\bigg(\psi_B - \bigg(\frac{\partial\psi_A}{\partial\theta}\bigg)^2 + \bigg(\frac{\partial\psi_A}{\partial\phi}\bigg)^2\bigg) + \opO(t^3).\eqnum{\nexteqnno[AsymptoticFormulaForF]}
$$
\endproclaim
\proclabel{AsymptoticFormulaForF}
In order to prove Lemma \procref{AsymptoticFormulaForF}, we first express $f$ in terms of the flows of certain Killing fields. Indeed, by the classical theory of symmetric spaces, the exponential map of $\opCL$ at $T_0$ is
$$
\Cal{E}:\opEnd(\Bbb{R}^2)\rightarrow\opCL;A\mapsto\opExp(\xi_A)T_0.\eqnum{\nexteqnno[ExponentialMapOfCL]}
$$
We may therefore suppose that
$$
\gamma(t) := \Cal{E}(tA+t^2B/2).\eqnum{\nexteqnno[FormulaForGamma]}
$$
Now let $X_A$, $X_B$, $F_A$ and $F_B$ be as in Section \subheadref{KillingFieldsOverCL}. Let $\alpha_{s,t}$ and $h_{s,t}$ be such that
$$
F_{B,s}(F_{A,t}(\theta,\phi,0)) = (\alpha_{s,t}(\theta,\phi),h_{s,t}(\theta,\phi)).\eqnum{\nexteqnno[DefinitionOfAlphaAndH]}
$$
Since $\alpha_{0,0}=\opId$, $\alpha_{s,t}$ is a smooth diffeomorphism for sufficiently small $(s,t)$. For all such $(s,t)$, let $\beta_{s,t}$ be its inverse.
\proclaim{Lemma \nextprocno}
\noindent With $f$ as above,
$$
f_t(\theta,\phi) = (h_{t^2/2,t}\circ\beta_{t^2/2,t})(\theta,\phi) + \opO(t^3).\eqnum{\nexteqnno[AsymptoticFormulaForFConstructingCurves]}
$$
\endproclaim
\proclabel{AsymptoticFormulaForFConstructingCurves}
\proof Indeed, for all $(s,t)$,
$$
\opExp(sA+tB) = \opExp(tB)\opExp(sA) + \opO(st),
$$
so that, for all sufficiently small $t$,
$$\eqalign{
\gamma(t) &= \opExp(t^2B/2)\opExp(tA)T_0 + \opO(t^3)\hfill\cr
&=\left\{ (\Phi\circ F_{B,t^2/2}\circ F_{A,t})(\theta,\phi,0) + \opO(t^3)\ |\ \theta,\phi\in S^1\right\}\hfill\cr
&=\left\{ \Phi(\alpha_{t^2/2,t}(\theta,\phi),h_{t^2/2,t}(\theta,\phi)) + \opO(t^3)\ |\ \theta,\phi\in S^1\right\}\hfill\cr
&=\left\{ \Phi(\theta,\phi,(h_{t^2/2,t}\circ\beta_{t^2/2,t})(\theta,\phi)) + \opO(t^3)\ |\ \theta,\phi\in S^1\right\}\hfill\cr
&=\ope[h_{t^2/2,t}\circ\beta_{t^2/2,t}] + \opO(t^3).\hfill\cr}
$$
The result follows.\qed
\medskip
\noindent Denote $\alpha_t:=\alpha_{0,t}$, $h_t:=h_{0,t}$ and $\beta_t:=\beta_{0,t}$.
\proclaim{Lemma \nextprocno}
\noindent The derivatives of $h$, $\alpha$ and $\beta$ satisfy
$$\eqalign{
\frac{\partial h_t}{\partial t}\bigg|_{t=0} &= \psi_A,\cr
\frac{\partial\alpha_t}{\partial t}\bigg|_{t=0} &= \bigg(\frac{\partial\psi_A}{\partial\theta}\bigg)\partial_\theta - \bigg(\frac{\partial\psi_A}{\partial\phi}\bigg)\partial_\phi,\cr
\frac{\partial \beta_t}{\partial t}\bigg|_{t=0} &= -\bigg(\frac{\partial\psi_A}{\partial\theta}\bigg)\partial_\theta + \bigg(\frac{\partial\psi_A}{\partial\phi}\bigg)\partial_\phi\ \text{and}\cr
\frac{\partial^2 h_t}{\partial t^2}\bigg|_{t=0} &= \bigg(\frac{\partial\psi_A}{\partial\theta}\bigg)^2 - \bigg(\frac{\partial\psi_A}{\partial\phi}\bigg)^2.\cr}\eqnum{\nexteqnno[DerivativesOfHAlphaBeta]}
$$
\endproclaim
\proclabel{DerivativesOfHAlphaBeta}
\proof Indeed, by definition, for all $t$, and for all $\theta,\phi$,
$$
(\alpha_t(\theta,\phi),h_t(\theta,\phi)) = F_{A,t}(\theta,\phi,0),
$$
so that
$$
\bigg(\frac{\partial\alpha_t}{\partial t},\frac{\partial h_t}{\partial t}\bigg) = X_A(\alpha_t,h_t).\eqnum{\nexteqnno[DerivativesOfHAlpha]}
$$
Since $\alpha_0=\opId$ and $h_0=0$, by \eqnref{ExplicitFormulaForX}, this yields
$$\eqalign{
\frac{\partial h_t}{\partial t}\bigg|_{s=0} &= \psi_A\ \text{and}\cr
\frac{\partial \alpha_t}{\partial t}\bigg|_{t=0} &= \bigg(\frac{\partial\psi_A}{\partial\theta}\bigg)\partial_\theta - \bigg(\frac{\partial\psi_A}{\partial\phi}\bigg)\partial_\phi.\cr}
$$
Furthermore, by the chain rule
$$
\frac{\partial\beta_t}{\partial t}\bigg|_{t=0} = -\frac{\partial\alpha_t}{\partial t}\bigg|_{t=0} = -\bigg(\frac{\partial\psi_A}{\partial\theta}\bigg)\partial_\theta+\bigg(\frac{\partial\psi_A}{\partial\phi}\bigg)\partial_\phi.
$$
Finally, by \eqnref{ExplicitFormulaForX} and \eqnref{DerivativesOfHAlpha},
$$
\frac{\partial h_t}{\partial t} = \psi_A\circ\alpha_t.
$$
Differentiating and evaluating at $0$ yields
$$
\frac{\partial^2 h_t}{\partial t^2}\bigg|_{t=0} = \bigg(\frac{\partial\psi_A}{\partial\theta}\bigg)^2 - \bigg(\frac{\partial\psi_A}{\partial\phi}\bigg)^2.
$$
This completes the proof.\qed
\medskip
\noindent We now prove Lemma \procref{AsymptoticFormulaForF}.
\medskip
{\bf\noindent Proof of Lemma \procref{AsymptoticFormulaForF}:\ }Indeed, by \eqnref{DerivativesOfHAlphaBeta},
$$\eqalign{
\frac{\partial h_{s,t}}{\partial s}\bigg|_{s,t=0} &= \psi_B\ \text{and}\cr
\frac{\partial h_{s,t}}{\partial t}\bigg|_{s,t=0} &= \psi_A.\cr}
$$
Thus, by \eqnref{AsymptoticFormulaForFConstructingCurves},
$$\eqalign{
\frac{\partial f_t}{\partial t}
&=t\bigg(\frac{\partial h}{\partial s}\bigg)_{t^2/2,t}\circ\beta_{t^2/2,t} + \bigg(\frac{\partial h}{\partial t}\bigg)_{t^2/2,t}\circ\beta_{t^2/2,t}\cr
&\qquad + \big(\opD h_{t^2/2,t}\circ\beta_{t^2/2,t}\big)\bigg(t\bigg(\frac{\partial\beta}{\partial s}\bigg)_{t^2/2,t} + \bigg(\frac{\partial\beta}{\partial t}\bigg)_{t^2/2,t}\bigg) + \opO(t^2).\cr}
$$
Since $h_{0,0}=0$ and $\beta_{0,0}=\opId$, evaluating at $t=0$ yields
$$
\frac{\partial f_t}{\partial t}\bigg|_{t=0} = \frac{\partial h}{\partial t}\bigg|_{s,t=0} = \psi_A.
$$
Differentiating a second time and evaluating at zero yields
$$
\frac{\partial^2 f_t}{\partial t^2}\bigg|_{t=0}
=\frac{\partial h}{\partial s}\bigg|_{s,t=0} + \frac{\partial^2 h}{\partial t^2}\bigg|_{s,t=0}
+2\opD\bigg(\frac{\partial h}{\partial t}\bigg)\frac{\partial\beta}{\partial t}\bigg|_{s,t=0},
$$
so that, by \eqnref{DerivativesOfHAlphaBeta} again
$$
\frac{\partial^2 f_t}{\partial t^2}\bigg|_{t=0} = \psi_B - \bigg(\frac{\partial\psi_A}{\partial\theta}\bigg)^2 + \bigg(\frac{\partial\psi_A}{\partial\phi}\bigg)^2.
$$
The result now follows by Taylor's theorem.\qed
\newsubhead{Zero sets of elements of $\opKer(\Delta+2)$}[ZeroSetsOfTheKernel]
We conclude this chapter by reviewing the geometries of the zero sets of elements of $\opKer(\Delta+2)$. We will study this in some detail since, not only is it useful in establishing the Morse-Smale property of functions constructed in the sequel, it also provides a charming geometric application of the Fermi parametrisation given in \eqnref{FermiParametrisation}.
\par
For all non-zero $A\in\opEnd(\Bbb{R}^2)$, denote
$$
Z_A := \psi_A^{-1}(\left\{0\right\}).\eqnum{\nexteqnno[DefinitionOfZA]}
$$
Observe that if
$$
J := \pmatrix 0\hfill& -1\hfill\cr 1\hfill& 0\hfill\cr\endpmatrix\eqnum{\nexteqnno[DefinitionOfJ]}
$$
denotes the standard complex structure of $\Bbb{R}^2$ then, for all $(\theta,\phi)$ and for all $(s,t)$,
$$
\psi_A(\theta+s,\phi+t) = \psi_{\opExp(-sJ)A\opExp(tJ)},\eqnum{\nexteqnno[RotationsYieldTranslations]}
$$
so that composing $A$ with rotations has the effect of rotating $Z_A$ in $\Bbb{T}$. In particular, $Z_A$ is invariant under half-turns of each of the two components of $\Bbb{T}:=\Bbb{S}^1\times\Bbb{S}^1$.
\par
Let $\Sigma$ now be the unit sphere in $\opEnd(\Bbb{R}^2)$ with respect to the metric \eqnref{MetricOverEnd}. It will trivially suffice to study only $A\in\Sigma$. Define the subsets $\Sigma_0$ and $\Sigma_{\pi/4}$ of $\Sigma$ by
$$\eqalign{
\Sigma_0 &:= \left\{ A\in\Sigma\ |\ \opDet(A)=0\right\}\ \text{and}\cr
\Sigma_{\pi/4} &:= \left\{ A\in\Sigma\ |\ AA^t = \opId/4\pi^2\right\}.\cr}\eqnum{\nexteqnno[DefinitionOfPortionsOfSigma]}
$$
Matrices in these subsets will be said to be {\sl singular} and {\sl special} respectively, whilst all other matrices in $\Sigma$ will be said to be {\sl generic}.
\par
Define $\Theta:\opEnd(\Bbb{R}^2)\rightarrow\Bbb{R}^4$ by
$$
\Theta(A) := (\Theta_+(A),\Theta_-(A)),\eqnum{\nexteqnno[DefinitionOfTheta]}
$$
where,
$$
\Theta_\pm\bigg(\pmatrix a\hfill& b\hfill\cr c\hfill& d\hfill\cr\endpmatrix\bigg):=\pi(a\pm d,b\mp c).\eqnum{\nexteqnno[DefinitionOfThetaPM]}
$$
This function is a linear isometry with respect to the metric \eqnref{MetricOverEnd} of $\opEnd(\Bbb{R}^2)$ and the standard metric over $\Bbb{R}^4$. In particular, $\Theta$ sends $\Sigma$ to the unit sphere, sends $\Sigma_0$ to the standard Clifford torus $T_0$, sends $\Sigma_{\pi/4}$ to the union of the two equatorial circles $\Bbb{S}^3\minter(\Bbb{R}^2\times\left\{(0,0)\right\})$ and  $\Bbb{S}^3\minter(\left\{(0,0)\right\}\times\Bbb{R}^2)$ and sends the space of diagonal matrices to the subspace
$$
D:=\left\{(x,0,y,0)\ |\ x,y\in\Bbb{R}\right\}.\eqnum{\nexteqnno[DefinitionOfD]}
$$
The intersection $\Bbb{S}^3\minter D$ is a great circle meeting $T_0$ orthogonally at precisely $4$ points. This reflects the fact that there are precisely $4$ singular diagonal matrices of unit norm, namely
$$
\frac{1}{\sqrt{2}\pi}\pmatrix \pm1\hfill& 0\hfill\cr 0\hfill& 0\hfill\cr\endpmatrix\ \text{and}\
\frac{1}{\sqrt{2}\pi}\pmatrix 0\hfill& 0\hfill\cr 0\hfill& \pm1\hfill\cr\endpmatrix.
$$
We now obtain the following decomposition result.
\proclaim{Lemma \nextprocno}
\noindent For all $A\in\opEnd(\Bbb{R}^2)$, there exists a diagonal matrix $\Delta$ and a pair $(M,N)$ of special orthogonal matrices such that
$$
A = M\Delta N.
$$
Furthermore, if $A$ is not a scalar multiple of an orthogonal matrix, then there exist precisely $4$ such $\Delta$, and for each such $\Delta$, the pair $(M,N)$ is unique up to sign.
\endproclaim
\proclabel{MatrixDecompositionLemma}
\proof Indeed, a long but straightforward calculation shows that, for all $M\in\opSO(2)$ and for all $A\in\opEnd(\Bbb{R}^2)$,
$$\eqalign{
\Theta(MAM^{-1}) &= (\Theta_+(A),M^2\Theta_-(A))\ \text{and}\cr
\Theta(MAM) &= (M^{-2}\Theta_+(A),\Theta_-(A)).\cr}\eqnum{\nexteqnno[EquivarianceOfTheta]}
$$
Now choose $A\in\opEnd(\Bbb{R}^2)$. If $A=\lambda M$ for some orthogonal matrix $M$, then existence follows trivially. Otherwise, upon multiplying by a scalar factor, we may suppose that $A$ is a non-special element of $\Sigma$. Consider now the actions $\sigma$ and $\rho$ of $\opSO(2)^2$ defined over $\Sigma$ and $\Bbb{S}^3$ respectively by
$$\eqalign{
\sigma(M,N)A &:= MAN\ \text{and}\cr
\rho(M,N)(x,y) &:= (MNx,MN^{-1}y).\cr}
$$
By \eqnref{EquivarianceOfTheta}, $\Theta$ intertwines $\sigma$ and $\rho$. However, the $\rho$-orbit of $x:=\Theta(A)$ corresponds to a horizontal slice in the Fermi parametrisation \eqnref{FermiParametrisation}. In particular, this orbit meets the great circle $\Bbb{S}^3\minter D$ in precisely $4$ points. Let $y:=\rho(M^{-1},N^{-1})(x)$ be one such point. Then $\Delta:=\Theta^{-1}(y)$ is diagonal and
$$
A = M\Delta N.
$$
Since the pair $(M,N)$ is uniquely defined up to sign, this completes the proof.\qed
\medskip
It is now straightforward to describe the geometry of $Z_A$. We consider singular and non-singular matrices separately.
\medskip
{\noindent\bf Singular matrices.\ }Choose $A\in\Sigma_0$. By Lemma \procref{MatrixDecompositionLemma}, we may suppose that
$$
A = \frac{1}{\sqrt{2}\pi}\pmatrix 1\hfill& 0\hfill\cr 0\hfill&0\hfill\cr\endpmatrix,\eqnum{\nexteqnno[StandardSingularMatrix]}
$$
so that
$$
\psi_A = \frac{1}{\sqrt{2}\pi}\opCos(\theta)\opCos(\phi),
$$
and therefore
$$
Z_A = \left\{\theta=\pi/2\right\}\munion\left\{\theta=3\pi/2\right\}\munion\left\{\phi=\pi/2\right\}\munion\left\{\phi=3\pi/2\right\}.
\eqnum{\nexteqnno[SingularZeroSet]}
$$
In other words $Z_A$ consists of the union of two horizontal and two vertical circles.
\medskip
{\noindent\bf Non-singular matrices.\ }Choose $A\in\Sigma\setminus\Sigma_0$. By Lemma \procref{MatrixDecompositionLemma}, we may suppose that
$$
A = \frac{1}{\sqrt{2}\pi}\pmatrix \opCos(r)\hfill& 0\hfill\cr 0\hfill& -\opSin(r)\hfill\cr\endpmatrix,\eqnum{\nexteqnno[StandardGenericMatrixII]}
$$
for some non-zero $r\in[-\pi/4,\pi/4]$. It follows that
$$
\psi_A(\theta,\phi) = \frac{1}{\sqrt{2}\pi}\opCos(r)\opCos(\theta)\opCos(\phi) - \frac{1}{\sqrt{2}\pi}\opSin(r)\opSin(\theta)\opSin(\phi),
$$
so that, upon substituting
$$
(x,y) := (\opCot(\theta),\opTan(\phi)),\eqnum{\nexteqnno[ZeroSetSubstitution]}
$$
we obtain
$$
Z_A\minter(]0,\pi[\times]-\pi/2,\pi/2[) = \left\{(x,y)\ |\ y=\opCot(r)x\right\}.\eqnum{\nexteqnno[NonSingularNonSpecialZeroSetII]}
$$
By \eqnref{RotationsYieldTranslations}, with the exception of $4$ points, $Z_A$ coincides with the union of $4$ translates of this set. In particular, $Z_A$ is the union of two smooth, closed curves which, up to the transformation \eqnref{ZeroSetSubstitution}, are straight lines. Finally, when $A\in\Sigma_{\pi/4}$ is special, $r=\pm\pi/4$, so that $Z_A$ in fact consists of the union of two diagonal circles with constant gradient equal to $\opSign(\opDet(A))$.
\medskip
Finally, we have
\proclaim{Lemma \nextprocno}
\myitem{(1)} If $A\neq 0$ is non-special then, for all $B$,
$$
\psi_B - \bigg(\frac{\partial\psi_A}{\partial\theta}\bigg)^2 + \bigg(\frac{\partial\psi_A}{\partial\phi}\bigg)^2\eqnum{\nexteqnno[SecondDerivative]}
$$
is non-vanishing over an open, dense subset of $Z_A$.
\medskip
\myitem{(2)} If $A\neq 0$ is special and if $B\neq 0$ is not colinear with $A$, then \eqnref{SecondDerivative} is non-vanishing over an open dense subset of $Z_A$.
\endproclaim
\proclabel{TheSecondDerivativeIsNonVanishing}
\proof By Lemma \procref{MatrixDecompositionLemma}, we may suppose that
$$\eqalign{
A &= \frac{1}{\sqrt{2}\pi}\pmatrix\opCos(r)\hfill& 0\hfill\cr 0\hfill& -\opSin(r)\hfill\cr\endpmatrix\ \text{and}\cr
B &= \pmatrix \alpha\hfill& \beta\hfill\cr \gamma\hfill& \delta\hfill\cr \endpmatrix.\cr}
$$
Substituting $x:=\opCot(\theta)$ and $y:=\opTan(\theta)$, over $Z_A$, we have
$$
y = \opCot(r)x,
$$
so that, over $Z_A$,
$$
\bigg(\frac{\partial\psi_A}{\partial\theta}\bigg)^2 - \bigg(\frac{\partial\psi_A}{\partial\phi}\bigg)^2
=\frac{1}{2\pi^2}\opSin^2(\theta)\opCos^2(\theta)\opCos(2r)(1 - \opCot^2(r)x^4),
$$
and
$$\eqalign{
\phi_B &= \opSin^2(\theta)\opCos^2(\theta)\sqrt{1+x^2}\sqrt{1+\opCot^2(r)x^2}\cr
&\qquad\qquad\times(\gamma + (\alpha+\delta\opCot(r))x + \beta\opCot(r)x^2).\cr}
$$
For $\left|r\right|<\pi/4$, \eqnref{SecondDerivative} is trivially non-vanishing over an open, dense subset of $Z_A$ and the first assertion follows. When $\left|r\right|=\pi/4$, \eqnref{SecondDerivative} is non-vanishing over an open, dense subset of $Z_A$ unless
$$
\beta = \gamma = \alpha + \delta\opCot(r) = 0.
$$
The second assertion follows. This completes the proof.\qed
\newhead{The Morse-Smale property}[TheMorseSmaleProperty]
\newsubhead{The Morse property}[ConstructingFunctionsOfMorseType]
Recall from the introduction that $\opI:C^\infty(\Bbb{S}^3)\rightarrow C^\infty(\opCL)$ is defined by
$$
\opI[u](T) := \int_Tu\opdArea_T,\eqnum{\nexteqnno[DefinitionOfI]}
$$
where, for all $T$, $\opdArea_T$ denotes its standard area form. In this chapter, we study the analytic properties of this function. We first show
\proclaim{Theorem \nextprocno}
\noindent There exists a generic subset $\Cal{U}_1$ of $C^\infty(\Bbb{S}^3)$ such that, for all $u\in\Cal{U}_1$, $\opI[u]$ is of Morse type.
\endproclaim
\proclabel{TheMorseProperty}
\noindent Define $\opJ_1:C^\infty(\Bbb{S}^3)\times\opCL\rightarrow \opT^*\opCL$ by
$$
\opJ_1[u](T) := d\opI[u](T).\eqnum{\nexteqnno[FirstJetOfU]}
$$
Observe that $\opI[u]$ is of Morse type if and only if $\opJ_1[u]$ is transverse to the zero section in $\opT^*\opCL$. Theorem \procref{TheMorseProperty} therefore follows from Lemma \procref{FirstSubmersionLemma}, below, together with the following version of the Sard-Smale theorem (c.f. \cite{SmiRos} and \cite{Smale}).
\proclaim{Lemma \nextprocno}
\noindent Let $\Omega$ be a weakly smooth manifold. Let $X$, $Y$ and $Z$ be finite-dimensional manifolds. Let $Z_0$ be a smooth submanifold of $Z$. Let $F:\Omega\times X\times Y\rightarrow Z$ be weakly smooth. If $D_1F\oplus D_3F$ is everywhere surjective, then there exists a generic subset $\Omega_0$ of $\Omega$ such that, for all $u\in\Omega_0$ and for all $x\in X$, $F(u,x,\cdot)$ is transverse to $Z_0$.
\endproclaim
\proclabel{SardSmaleTheorem}
{\bf\noindent Sketch of proof:\ } We may suppose that $\Omega=C^\infty(W)$ for some compact manifold $W$. Let $X'$ and $Y'$ be relatively compact open subsets of $X$ and $Y$ respectively. Let $\Omega'$ be the set of all $u\in\Omega$ such that, for all $x\in X'$, the restriction of $F(u,x,\cdot)$ to $Y'$ is transverse to $Z_0$. By compactness of $X'$ and $Y'$, $\Omega'$ is open. We now show that this subset is dense. Choose $u\in C^\infty(W)$. Given any vector subspace $E$ of $C^\infty(W)$ and any $\epsilon>0$, denote
$$
E_{u,\epsilon} := B_\epsilon(u)\minter(u+E),
$$
where $B_\epsilon(u)$ here denotes the ball of radius $\epsilon$ about $u$ with respect to some norm. By compactness, there exists a finite dimensional vector subspace $E$ of $C^\infty(W)$ and $\epsilon>0$ such that the restriction of $D_1F\oplus D_3F$ to $E\oplus T_yY$ is surjective at every point of $E_{u,\epsilon}\times X'\times Y'$. It now follows by classical differential topology that there exists a generic subset $\Cal{U}$ of $E_{u,\epsilon}$ such that, for all $u'\in\Cal{U}$ and for all $x\in X'$, the restriction of $F(u,x,\cdot)$ to $Y'$ is transverse to $Z_0$. Since $u\in\Omega$ is arbitrary, it follows that $\Omega'$ is dense, as asserted. The result now follows upon taking compact exhaustions of $X$ and $Y$.\qed
\medskip
In much of what follows, the surjectivity results required to apply Lemma \procref{SardSmaleTheorem} will be derived from the following
\proclaim{Lemma \nextprocno}
\noindent For all $T\in\opCL$, for all $k\geq 0$, for all $\psi\in\opT_T\opCL\setminus\left\{0\right\}$ and for all $\alpha\in\opT^*_T\opCL$, there exists $u\in C^\infty(\Bbb{S}^3)$ such that, for all $0\leq i\leq k-1$,
$$
\nabla^i d\opI[u](T)(\psi,...,\psi,\cdot) = 0,
$$
and
$$
\nabla^k d\opI[u](T)(\psi,...,\psi,\cdot) = \alpha.
$$
\endproclaim
\proclabel{PreSurjectivity}
\remark Observe that when $k=0$, the hypothesis on $\psi$ is unnecessary.
\medskip
\proof Upon applying an element of $\opO(4)$, we may suppose that $T=T_0$. Let $A\in\opEnd(\Bbb{R}^2)$ be such that $\psi=\psi_A$. Let $\Omega\subseteq\Bbb{T}$ be an open subset over which $\psi_A$ is strictly positive. Let $g\in C_0^\infty(\Omega)$ be such that, for all $B\in\opEnd(\Bbb{R}^2)$,
$$
\alpha(\psi_B) = \frac{1}{2}\int_{\Bbb{T}}g\psi_Bd\theta d\phi.
$$
Let $\chi\in C^\infty_0(]-\pi/4,\pi/4[)$ be equal to $1$ near $0$. Define $\tilde{u}:\Bbb{T}\times]-\pi/4,\pi/4[\rightarrow\Bbb{R}$ by
$$
\tilde{u}(\theta,\phi,t) := \frac{t^kg(\theta,\phi)\chi(t)}{k!\psi_A(\theta,\phi)^{k-1}}.
$$
Observe that $\tilde{u}$ is smooth with compact support. Define $u$ over $\Bbb{S}^3$ by
$$
u(x) := \left\{
\matrix (\tilde{u}\circ\Phi^{-1})(x)\hfill&\ \text{if}\ x\in\opIm(\Phi)\ \text{and}\hfill\cr
0\hfill&\ \text{if}\ x\in\Phi(\opSupp(\tilde{u}))^c.\hfill\cr\endmatrix\right.
$$
We claim that $u$ is the desired function. Indeed, by the classical theory of symmetric spaces, the exponential map of $\opCL$ at $T_0$ is
$$
\Cal{E}:\opEnd(\Bbb{R}^2)\rightarrow\opCL;M\mapsto\opEnd(\xi_M)T_0.
$$
Choose $B\in\opEnd(\Bbb{R}^2)$ and define $\gamma:]-\epsilon,\epsilon[^k\rightarrow\opCL$ by
$$
\gamma(t_1,...,t_k) := \Cal{E}(t_1A + ... + t_{k-1}A + t_kB).
$$
Observe that all derivatives of $\opI[u]\circ\gamma$ up to and including order $(k-1)$ vanish at $0$. For all $(t_1,...,t_k)$, let $f_{t_1,...,t_k}$ be such that
$$
\gamma_{t_1,...,t_k} = e[f_{t_1,...,t_k}].
$$
Then, by \eqnref{AsymptoticFormulaForF} and Taylor's theorem,
$$
f_{t_1,...,t_k} = t_1\psi_A + ... + t_{k-1}\psi_A + t_k\psi_B + \opO((\left|t_1\right|+...+\left|t_k\right|)^2).
$$
For all $t_1,...,t_k$, let $\opdA_{t_1,...,t_k}$ be the area form of the graph of $f_{t_1,...,t_k}$. Then,
$$
\opdA_{t_1,...,t_k} = \frac{1}{2}d\theta d\phi + \opO(\left|t_1\right|+...+\left|t_k\right|).
$$
Combining these relations yields
$$\multiline{
\nabla^{k-1} d\opI[u](\psi_A,...,\psi_A,\psi_B)\cr
\qquad\qquad=\frac{\partial^k}{\partial t_1...\partial t_k}\opI[u](\gamma_{t_1,...,t_k})\bigg|_{t_1,...,t_k=0}\cr
\qquad\qquad=\frac{\partial^k}{\partial t_1...\partial t_k}\int_{\Bbb{T}}(u\circ\Phi)(\theta,\phi,f_{t_1,...,t_k}(\theta,\phi))\opdA_{t_1,...,t_k}\bigg|_{t_1,...,t_k=0}\cr
\qquad\qquad=\int_{\Bbb{T}}\frac{1}{2}g(\theta,\phi)\psi_B(\theta,\phi)d\theta d\phi\cr
\qquad\qquad=\beta(\psi_B).\vphantom{\frac{1}{2}}\cr}
$$
Since $B\in\opEnd(\Bbb{R}^2)$ is arbitrary, it follows that $u$ is indeed the desired function. This completes the proof.\qed
\medskip
\noindent In the present case, this yields
\proclaim{Lemma \nextprocno}
\noindent $\opJ_1$ is a submersion.
\endproclaim
\proclabel{FirstSubmersionLemma}
\proof Choose $(u,T)\in C^\infty(\Bbb{S}^3)\times\opCL$. Let $\opV$ be the vertical subbundle of $\opT\opT^*\opCL$. Let $\opD\opJ_1$ be the total derivative of $\opJ_1$ and, for each $i$, let $\opD_i\opJ_1$ be its partial derivative with respect to the $i$'th component. Since $\opJ_1[u]$ is a section of $\opT^*\opCL$,
$$
\big(\opD_2\opJ_1\cdot \opT_T\opCL\big)\oplus\opV = \opT\opT^*\opCL.
$$
By Lemma \procref{PreSurjectivity}, $\opV$ is also in the image of $\opD\opJ_1$, so that $\opD\opJ_1$ is surjective at this point. Since $(u,T)\in C^\infty(\Bbb{S}^3)\times\opCL$ is arbitrary, it follows that $\opJ_1$ is a submersion. This completes the proof.\qed
\medskip
\noindent This yields Theorem \procref{TheMorseProperty}.
\medskip
{\bf\noindent Proof of Theorem \procref{TheMorseProperty}:\ }Indeed, since $\opI[u]$ is of Morse type if and only if $\opJ_1[u]$ is transverse to the zero section, the result follows by Lemma \procref{FirstSubmersionLemma} and Theorem \procref{SardSmaleTheorem} upon setting $\Omega:=C^\infty(\Bbb{S}^3)$, $X:=\left\{0\right\}$ is the manifold consisting of a single point, $Y:=\opCL$, $Z:=T^*\opCL$, and $Z_0$ is the graph of the zero section in $\opT^*\opCL$.\qed
\newsubhead{The Morse-Smale property}[TheSmaleProperty]
The remainder of this chapter is devoted to proving the Morse-Smale property of $\opI[u]$ for generic $u$. Recall from the introduction that a smooth function over $\opCL$ is said to be of {\sl Morse-Smale type} whenever, in addition to it being of Morse type, every unstable manifold of its gradient flow is transverse to every stable manifold of this flow. We now describe an alternative, more technical, characterisation of the Morse-Smale property which is better adapted to our current applications. We refer the reader to \cite{Schwarz} for a complete treatment of the relevant theory.
\par
First, we define the {\sl gradient flow operator} $\opGF:C^1(\Bbb{R},\opCL)\times C^1(\opCL)\rightarrow C^0(\Bbb{R},T\opCL)$ by
$$
\opGF[\gamma,f](t) := \dot{\gamma}(t) + \nabla f(\gamma(t)).\eqnum{\nexteqnno[GradientFlowOperator]}
$$
Given a function $f\in C^1(\opCL)$, a smooth curve $\gamma:\Bbb{R}\rightarrow\opCL$ is a complete {\sl gradient flow} of $f$ whenever $\opGF[\gamma,f]$ vanishes.
\par
Consider a point $(\gamma,f)\in C^1(\Bbb{R},\opCL)\times C^1(\opCL)$. The partial derivative of $\opGF$ with respect to the first component at this point defines a first-order, linear differential operator $\opD_1\opGF[\gamma,f]$ from $\Gamma^1(\gamma^*T\opCL)$ into $\Gamma^0(\gamma^*T\opCL)$. We obtain an explicit formula for this operator as follows. Denote $E:=\opT_{\gamma(0)}\opCL$. Let $\tau:\Bbb{R}\times E\rightarrow\gamma^*\opT\opCL$ be parallel transport along $\gamma$. Define the matrix-valued function $M_{\gamma,f}:\Bbb{R}\rightarrow\opEnd(E)$ such that, for all $t\in\Bbb{R}$ and for all $\xi,\eta\in E$,
$$
\langle M_{\gamma,f}(t)\xi,\eta\rangle := \opHess(f)(\gamma(t))(\tau(t)\xi,\tau(t)\eta).\eqnum{\nexteqnno[DefinitionOfM]}
$$
Elements of $C^1(\Bbb{R},E)$ identify through $\tau$ with elements of $\Gamma^1(\gamma^*\opT\opCL)$ which in turn identify with first order variations of $\gamma$. The action of $\opD_1\opGF[\gamma,f]$ on such first-order variations yields elements of $\Gamma^0(\gamma^*\opT\opCL)$ which in turn identify through $\tau^{-1}$ with elements of $C^0(\Bbb{R},E)$. Following through these identifications, we obtain, for all $g\in C^1(\Bbb{R},E)$,
$$
(\opD_1\opGF[\gamma,f]g)(t) = \dot{g}(t) - M_{\gamma,f}(t)g(t),\eqnum{\nexteqnno[DOneGradientFlowOperator]}
$$
\par
When $f$ is of Morse type, as $t$ tends to $\pm\infty$, $\gamma(t)$ converges to critical points of $\nabla f$ and $M_{\gamma,f}(t)$ converges to non-singular matrices $M_{\gamma,f}(\pm\infty)$. By the Atiyah-Patodi-Singer index theorem (see \cite{Salamon}), for all $\alpha\in]0,1[$, $\opD_1\opGF[\gamma,f]$ defines a Fredholm map from $C^{1,\alpha}(\Bbb{R},E)$ into $C^{0,\alpha}(\Bbb{R},E)$ of Fredholm index equal to the Morse index of $M_{\gamma,f}(+\infty)$ minus that of $M_{\gamma,f}(-\infty)$. The function $f$ is now said to be of {\sl Morse-Smale type} whenever, in addition to it being of Morse type, the operator $\opD_1\opGF[\gamma,f]$ is surjective for every gradient flow $\gamma$ of $f$. In the rest of this chapter, we show
\proclaim{Theorem \nextprocno}
\noindent There exists a generic subset $\Cal{U}_2$ of $C^\infty(\Bbb{S}^3)$ such that, for all $u\in\Cal{U}_2$, $\opI[u]$ is of Morse-Smale type.
\endproclaim
\proclabel{TheMorseSmaleProperty}
As with Theorem \procref{TheMorseProperty}, Theorem \procref{TheMorseSmaleProperty} follows from a suitable surjectivity result together with a suitable version of the Sard-Smale theorem. The version of the Sard-Smale theorem used is more technical than Lemma \procref{SardSmaleTheorem}, not only because $X$, $Y$ and $Z$ are now infinite-dimensional Banach manifolds, but also because the precompactness properties of families of gradient flows are actually quite subtle. In what follows, we will satisfy ourselves with the surjectivity result, and we refer the reader again to \cite{Schwarz} for a complete treatment of the remaining part of the theory.
\par
Given a point $(\gamma,f)\in C^1(\Bbb{R},\opCL)\times C^1(\opCL)$, the partial derivative of $\opGF$ with respect to the second component at this point defines a first-order, linear, partial differential operator $\opD_2\opGF[\gamma,f]$ from $C^1(\opCL)$ to $\Gamma^0(\gamma^*T\opCL)$. With the above identifications, we have, for all $g\in C^1(\opCL)$,
$$
(\opD_2\opGF[\gamma,f]g)(t) = \tau(t)^{-1}\nabla g(\gamma(t)).\eqnum{\nexteqnno[DTwoGradentFlowOperator]}
$$
Bearing in mind the discussion of the preceding paragraph, Theorem \procref{TheMorseSmaleProperty} follows from
\proclaim{Lemma \nextprocno}
\noindent There exists a generic subset $\Cal{U}_3$ of $C^\infty(\Bbb{S}^3)$ such that, for all $u\in\Cal{U}_3$ and for every gradient flow $\gamma$ of $\opI[u]$, the sum
$$
\opD_1\opGF[\gamma,\opI[u]]+\opD_2\opGF[\gamma,\opI[u]]\opD\opI[u]
$$
defines a surjective linear map from $C^{1,\alpha}(\Bbb{R},E)\oplus C^\infty(\Bbb{S}^3)$ into $C^{0,\alpha}(\Bbb{R},E)$.
\endproclaim
\proclabel{SurjectivityMorseSmale}
Finally, we apply Fredholm theory to transform Lemma \procref{SurjectivityMorseSmale} into the form that we will prove.
Choose $u\in\Cal{U}_1$ so that $\opI[u]$ is of Morse type. Let $\gamma:\Bbb{R}\rightarrow\opCL$ be a gradient flow of $\opI[u]$. Since $M_{\gamma,\opI[u]}(t)$ is symmetric for all $t$, the formal dual of $\opD_1\opGF[\gamma,\opI[u]]$ satisfies
$$
(\opD_1\opGF[\gamma,\opI[u]]^*g)(t) := -\dot{g}(t) - M_{\gamma,\opI[u]}(t)g(t).\eqnum{\nexteqnno[DualDOneGradientFlowOperator]}
$$
Its kernel is a finite-dimensional space of smooth functions which decay exponentially at infinity (see \cite{Schwarz}). In particular,
$$
\opKer(\opD_1\opGF[\gamma,\opI[u]]^*) \subseteq L^1(\Bbb{R},E).
$$
Since $\opD_2\opGF[\gamma,\opI[u]]g\subseteq L^\infty(\Bbb{R},E)$ for all $g$, it follows that the pairing
$$
\opKer(\opD_1\opGF[\gamma,\opI[u]]^*)\oplus C^\infty(\opCL)\rightarrow\Bbb{R};
(g,h)\mapsto\int_{-\infty}^{\infty}\langle g,D_2\opGF[\gamma,\opI[u]]h\rangle dt\eqnum{\nexteqnno[PairingToProveSurjectivity]}
$$
is well-defined. Lemma \procref{SurjectivityMorseSmale} is thus equivalent to
\proclaim{Lemma \nextprocno}
\noindent There exists a generic subset $\Cal{U}_3$ of $C^\infty(\Bbb{S}^3)$ such that, for all $u\in\Cal{U}_3$, for every gradient flow $\gamma$ of $\opI[u]$, and for all $g\in\opKer(\opD_1\opGF[\gamma,\opI[u]]^*)$, there exists $v\in C^\infty(\Bbb{S}^3)$ such that
$$
\langle g,D\opI[u]v\rangle \neq 0.
$$
\endproclaim
\proclabel{SurjectivityMorseSmaleII}
\noindent Lemma \procref{SurjectivityMorseSmaleII} will follow from Lemmas \procref{ConstructingVPartI}, \procref{ConstructingVPartII} and \procref{ConstructingVPartIII} of Section \subheadref{ConstructingThePerturbation}, below. Although the construction of $v$, carried out in Section \subheadref{ConstructingThePerturbation}, below, is relatively straightforward, showing that this function has the desired properties will require very careful estimates that will occupy most of the rest of this chapter.
\par
We conclude this section by establishing the following key property of function which lie in the kernel of $\opD_1\opGF[\gamma,\opI[u]]^*$.
\proclaim{Lemma \nextprocno}
\noindent For all $g\in\opKer(\opD_1\opGF[\gamma,\opI[u]]^*)$ and for all $t\in\Bbb{R}$, $g(t)$ is orthogonal to $\tau(t)^{-1}\dot{\gamma}(t)$.
\endproclaim
\proclabel{GIsOrthogonalToFlow}
\proof Indeed, let $\phi$ and $\tilde{g}$ be such that, for all $t$,
$$\eqalign{
g(t) &= \phi(t)\tau(t)^{-1}\dot{\gamma}(t) + \tilde{g}(t)\ \text{and}\cr
\tilde{g}(t) &\perp \tau(t)^{-1}\dot{\gamma}(t).\cr}
$$
The function $\phi$ and all its derivatives decay exponentially as $t$ tends to $\pm\infty$. For small $s$, define $\gamma_s:\Bbb{R}\rightarrow\opCL$ by
$$
\gamma_s(t) := \gamma\bigg(t+s\int_{-\infty}^t\phi(r) dr\bigg).
$$
The family $(\gamma_s)_{s\in]-\epsilon,\epsilon[}$ defines a smooth curve in $C^{1,\alpha}(\Bbb{R},\opCL)$ such that, for all $s$,
$$
\opGF[\gamma_s,\opI[u]](t) = s\phi(t)\dot{\gamma}\bigg(t+s\int_{-\infty}^t\phi(r)dr\bigg).
$$
Consequently
$$
\bigg(\opD_1\opGF[\gamma,\opI[u]]\frac{\partial\gamma_s}{\partial s}\bigg|_{s=0}\bigg)(t) = \phi(t)\tau(t)^{-1}\dot{\gamma}(t),
$$
so that
$$\eqalign{
\int_{-\infty}^\infty\left|\phi(t)\right|^2\|\dot{\gamma}(t)\|^2dt
&=\int_{-\infty}^\infty\bigg\langle\bigg(\opD_1\opGF[\gamma,\opI[u]]\frac{\partial\gamma_s}{\partial s}\bigg|_{s=0}\bigg)(t),f(t)\bigg\rangle dt\cr
&=\int_{-\infty}^\infty\bigg\langle\frac{\partial\gamma_s}{\partial s}\bigg|_{s=0}(t),\big(\opD_1\opGF[\gamma,\opI[u]]^*f\big)(t)\bigg\rangle dt\cr
&=\vphantom{\frac{1}{2}}0.\cr}
$$
It follows that $\phi=0$. Since $g\in\opKer(\opD_1\opGF[\gamma,\opI[u]])$ is arbitrary, this completes the proof.\qed
\newsubhead{Towards the Morse-Smale property}[TowardsTheMorseSmaleProperty]
The main challenge in constructing the function $v$ in Lemma \procref{SurjectivityMorseSmaleII} arises from the fact that any two Clifford tori will have non-trivial intersection in $\Bbb{S}^3$. Consequently, given any smooth curve $\gamma:\Bbb{R}\rightarrow\opCL$, constructing $v$ so as to have some desired effect near a given point of this curve may well, in fact, produce a self-cancelling effect near some other point. The following technical results serve to ensure that, along the curves of interest to us, such intersections are as infrequent and as regular as possible. For $k\geq 2$, define $\opJ_k:C^\infty(\Bbb{S}^3)\times\opCL\rightarrow\opT^*\opCL^{\oplus k}$ by
$$
\opJ_k[u](T) := (d\opI[u](T),...,\nabla^{k-1}d\opI[u](T)(d\opI[u](T),...,d\opI[u](T),\cdot)).\eqnum{\nexteqnno[HigherOrderJets]}
$$
\proclaim{Lemma \nextprocno}
\noindent For all $k\geq 2$, $\opJ_k$ is a submersion away from $J_1^{-1}(0)$, where $0$ here denotes the zero section in $T^*\opCL$.
\endproclaim
\proclabel{SecondSubmersionLemma}
\proof Choose $(u,T)\in C^\infty(\Bbb{S}^3)\times\opCL$ such that $J_1[u](T)\neq 0$. Let $\opV$ be the vertical subbundle of $\opT(\opT^*\opCL^{\oplus k})$. Observe that $\opV$ naturally decomposes as
$$
\opV = \opV_1\oplus...\oplus\opV_k,
$$
where, for all $(\alpha_1,...,\alpha_k)\in\opT^*_T\opCL^{\oplus k}$, the fibres of $\opV_1$,...,$\opV_k$ over this point each naturally identify with $\opT^*_T\opCL$. Since $\opJ_k[u]$ is a section of $\opT^*\opCL^{\oplus k}$, we have
$$
\big(D_2\opJ_k\cdot T_T\opCL\big)\oplus\opV_1\oplus...\oplus\opV_k = T(T^*\opCL^{\oplus k}).
$$
However, by Lemma \procref{PreSurjectivity}, $\opV_1\oplus...\oplus\opV_k$ is also in the image of $\opD\opJ_k$ so that $\opD\opJ_k$ is surjective at this point. Since $(u,T)\notin\opJ_1^{-1}(0)$ is arbitrary, it follows that $\opJ_k$ is a submersion away from this set. This completes the proof.\qed
\medskip
Recall (c.f. Section \subheadref{TheGeometryOfCurvesInCL}) that the space of rank $1$ matrices $A$ in $\opEnd(\Bbb{R}^2)$ is a $3$-dimensional submanifold diffeomorphic to $\Bbb{T}\times]0,\infty[$. For all such $A$, the zero set $\opZ_A:=\psi_A^{-1}(\left\{0\right\})$ consists of the union of two horizontal and two vertical circles in $\Bbb{T}$ whose intersections define $4$ points of $\Bbb{T}$. Let $\opZ_A^\opreg$ be the complement of these $4$ points in $\opZ_A$. For all other non-zero $A$, $\opZ_A$ is a union of finitely many smooth, closed curves in $A$ and we thus denote $\opZ_A^\opreg:=\opZ_A$.
\par
For every point $x\in T_0$, the set of tori in $\opCL$ which are tangent to $T_0$ at $x$ is a closed $1$-dimensional submanifold. Furthermore, for any such torus $T$, the points of tangency with $T_0$ are isolated unless $T=T_0$. It follows that, for all non-zero $A\in\opEnd(\Bbb{R}^2)$, the set $\opW^1_A$ of tori in $\opCL$, distinct from $T_0$, which intersect $T_0$ tangentially at some point of $\opZ_A^\opreg$ is a non-compact, codimension $2$ submanifold of $\opCL$. For all $A$, let $\opT\opW^1_A$ be the total space of the tangent bundle of $\opW^1_A$ in $\opT\opCL$. This is a non-compact, codimension $4$ submanifold of $\opT\opCL$. Lemmas \procref{SardSmaleTheorem} and \procref{SecondSubmersionLemma} now yield
\proclaim{Lemma \nextprocno}
\noindent There exists a generic subset $\Cal{U}_4$ of $C^\infty(\Bbb{S}^3)$ such that, for all $u\in\Cal{U}_4$, for all non-zero $A\in\opEnd(\Bbb{R}^2)$ and for all $M\in\opO(4)$, the graph of $\opJ_1[u]$ is transverse to $M\opT\opW^1_A$, where $M\opT\opW^1_A$ here denotes the image of $\opT\opW^1_A$ under $M$.
\endproclaim
\proclabel{FirstTransversalProperty}
When $A=\lambda M$ for some $M\in\opO(2)$ and some non-zero scalar $\lambda$, the set $\opZ_A$ consists of two diagonal lines with unit slope. Recall that the vector $\psi_A$ is said to be {\sl special} in this case. Define the submanifold $\opSpec\subseteq\opT^*\opCL\oplus\opT^*\opCL$ by
$$
\opSpec := \left\{ (T,\psi^\flat,\lambda\psi^\flat)\ |\ T\in\opCL,\ \psi\ \text{special},\ \lambda\in\Bbb{R}\right\},
$$
where the superscript $\flat$ here denotes Berger's musical isomorphism. Since the space of all special matrices is $2$-dimensional, it follows that, away from the zero section, $\opSpec$ is a codimension $5$ submanifold of $T^*\opCL\oplus T^*\opCL$. Lemmas \procref{SardSmaleTheorem} and \procref{SecondSubmersionLemma} now yield
\proclaim{Lemma \nextprocno}
\noindent There exists a generic subset $\Cal{U}_5$ of $C^\infty(\Bbb{S}^3)$ such that, for all $u\in\Cal{U}_5$, the graph of $\opJ_2[u]$ only intersects $\opSpec$ at zeroes of $\opJ_1[u]$.
\endproclaim
\proclabel{SecondTransversalProperty}
Finally, for every point $x\in T_0$, the set $\opW^2_x$ of tori in $\opCL$ which intersect $T_0$ transversally at $x$ is a non-compact, codimension $1$ submanifold of $\opCL$. Define $\opT^4\opW^2_x\subseteq\opT\opCL^{\oplus 4}$ by
$$
\opT^4\opW^2_x := \left\{(\gamma(0),\dot{\gamma}(0),...,(\nabla_{\dot{\gamma}}^3\dot{\gamma})(0))\ |\ \gamma:]-\epsilon,\epsilon[\rightarrow\opW^2_x\right\}.
$$
Observe that $\opT^4\opW^2_x$ is a smooth, codimension $5$ submanifold of $\opT\opCL^{\oplus 4}$. Lemmas \procref{SardSmaleTheorem} and \procref{SecondSubmersionLemma} now yield
\proclaim{Lemma \nextprocno}
\noindent There exists a generic subset $\Cal{U}_6$ of $C^\infty(\Bbb{S}^3)$ such that, for all $u\in\Cal{U}_6$, for all $x\in T_0$ and for all $M\in\opO(4)$, the graph of $\opJ_4[u]$ only intersects $M\opT^4\opW^2_x$ at zeroes of $\opJ_1[u]$, where $M\opT^4\opW^2_x$ here denotes the image of $\opT^4\opW^2_x$ under $M$.
\endproclaim
\proclabel{ThirdTransversalProperty}
\newsubhead{Constructing $v$}[ConstructingThePerturbation]
We now construct the function $v$ of Lemma \procref{SurjectivityMorseSmaleII}. We continue to use the notation of the preceding sections. Choose
$$
u\in\Cal{U}_3:=\Cal{U}_1\minter\Cal{U}_4\minter\Cal{U}_5\minter\Cal{U}_6,\eqnum{\nexteqnno[DefinitionOfUThree]}
$$
where $\Cal{U}_1$, $\Cal{U}_4$, $\Cal{U}_5$ and $\Cal{U}_6$ are as in Theorem \procref{TheMorseProperty} and Lemmas \procref{FirstTransversalProperty}, \procref{SecondTransversalProperty} and \procref{ThirdTransversalProperty} respectively. Let $\gamma:\Bbb{R}\rightarrow\opCL$ be a gradient flow of $\opI[u]$. Upon applying an element of $\opO(4)$, we may suppose that $\gamma(0)=T_0$. Let $A,B\in\opEnd(\Bbb{R}^2)$ be such that
$$\eqalign{
\dot{\gamma}(0) &= \psi_A\ \text{and}\cr
\big(\nabla_{\dot{\gamma}}\dot{\gamma}\big)(0) &= \psi_B.\cr}
$$
For all $t$, let $N_t$ be the unit normal vector field over $\gamma(t)$, let $J_t$ be its Jacobi operator and let $g_t\in\opKer(J_t)$ be the element corresponding to the vector $g(t)$. For sufficiently small $\epsilon$, let $(f_t)_{t\in]-\epsilon,\epsilon[}$ be such that, for all $t$,
$$
e[f_t] = \gamma(t).
$$
\par
Since $u\in\Cal{U}_1$, $\opI[u]$ is of Morse type. Since $u\in\Cal{U}_4$, $\gamma(t)$ only intersects $\opW^1_A$ over a discrete subset of $\Bbb{R}\setminus\left\{0\right\}$. In particular, the set of points of $\opZ^\opreg_A$ at which $\gamma(t)$ intersects $T_0$ tangentially for some $t$ is countable. There therefore exists a point $x_0\in\opZ^\opreg_A$ at which, for all $t\in\Bbb{R}\setminus\left\{0\right\}$, $\gamma(t)$ only intersects $T_0$ transversally. We may suppose, furthermore, that the limit tori $\gamma(\pm\infty)$ also only intersect $T_0$ transversally at this point. Since, by Lemma \procref{GIsOrthogonalToFlow}, $g(0)$ is orthogonal to $\psi_A$, upon modifying $x_0$ if necessary, we may suppose that
$$
g_0 := g(x_0) > 0.\eqnum{\nexteqnno[PositivityOfG]}
$$
Since $u\in\Cal{U}_5$, either $\psi_A$ is non-special or $\psi_A$ and $\psi_B$ are non-colinear. By Lemma \procref{TheSecondDerivativeIsNonVanishing}, upon modifying $x_0$ again if necessary, we may suppose that
$$
\frac{\partial^2 f_t}{\partial t^2}(x_0)\bigg|_{t=0} \neq 0.\eqnum{\nexteqnno[NonVanishingOfSecondDerivative]}
$$
Finally, upon applying another element of $\opO(4)$, we may suppose that $x_0=(0,0)$ and, since all calculations in the sequel are local in a neighbourhood of this point, we may suppose that the tangent line to $\opZ^\opreg_A$ at $(0,0)$ coincides with the $\theta$-axis.
\par
Choose $0<\delta\ll\eta\ll 1$. The parameter $\eta$ will be fixed presently and $\delta$ will be made to tend to $0$. For $r>0$, define
$$
\chi^1_r(x) :=
\left\{
\matrix 1\hfill&\ \text{if}\ \left|x\right|\leq r,\hfill\cr
2-\left|x\right|/r\hfill&\ \text{if}\ r<\left|x\right|\leq 2r\ \text{and}\hfill\cr
0\hfill&\ \text{otherwise,}\hfill\cr\endmatrix\right.\eqnum{\nexteqnno[FirstComponentOfV]}
$$
and
$$
\chi^2_r(x) :=
\left\{
\matrix x/r\hfill&\ \text{if}\ 0<x\leq r,\hfill\cr
2-x/r\hfill&\ \text{if}\ r<x\leq 2r\ \text{and}\hfill\cr
0\hfill&\ \text{otherwise.}\hfill\cr\endmatrix\right.\eqnum{\nexteqnno[SecondComponentOfV]}
$$
Define
$$
\tilde{v}_\delta := \tilde{v}_\delta(\theta,\phi,t) := \frac{1}{\delta}\chi^1_\delta(\theta)\chi^1_\delta(\phi)\chi^2_{\delta^3}(t),\eqnum{\nexteqnno[DefinitionOfTildeVDelta]}
$$
and
$$
v_\delta(x)
:=\left\{
\matrix (v\circ\Phi^{-1})(x)\hfill&\ \text{if}\ x\in\opIm(\Phi)\ \text{and}\hfill\cr
0\hfill&\ \text{if}\ x\in\Phi(\opSupp(\tilde{v}_\delta))^c.\hfill\cr\endmatrix\right.\eqnum{\nexteqnno[DefinitionOfVDelta]}
$$
We now show that, for sufficiently small $\delta$, every function sufficiently close to $v_\delta$ in $W^{1,1}(\Bbb{S}^3)$ has the desired properties.
\proclaim{Lemma \nextprocno}
\noindent For all $\eta>0$, there exists $C_1>0$ such that, for all $\delta<\eta$ and for all $\left|t\right|>\eta$,
$$
\|\nabla\opI[v_\delta](\gamma(t))\|\leq C_1.\eqnum{\nexteqnno[BoundsAwayFromZero]}
$$
\endproclaim
\proclabel{BoundsAwayFromZero}
\proof By definition of $x$, for all $t\in[-\infty,\infty]\setminus]-\eta,\eta[$, $\gamma(t)$ only intersects $\gamma(0)$ transversally at $x$. By compactness, there exists $\theta_0>0$ such that, if $\gamma(t)$ intersects $\gamma(0)$ at $x$, then the tangent planes of these two surfaces make an angle of at least $\theta_0$ with one another at this point. There therefore exists $A_1>0$ such that, for all $\delta<\eta$ and for all $t\in[-\infty,\infty]\setminus]-\eta,\eta[$,
$$
\Cal{H}_2(\gamma(t)\minter\opSupp(v_\delta))\leq A_1\delta^4,
$$
where $\Cal{H}_2$ here denotes $2$-dimensional Hausdorff measure. It follows that there exists $A_2>0$ such that, for all $\delta<\eta$, for all $t\in[-\infty,\infty]\setminus]-\eta,\eta[$ and for all $g\in C^\infty(\gamma(t))$,
$$
\int_{\gamma(t)}gdv_\delta(N_t)\opdA \leq A_2\|g\|_{L^\infty}.
$$
The result follows.\qed
\proclaim{Lemma \nextprocno}
\noindent For all $\eta>0$ and for all $\epsilon>0$, there exists $R>0$ such that, for all $\delta<\eta$,
$$
\left|\int_{\Bbb{R}\setminus[-R,R]}\langle g(t),\tau(t)^{-1}\nabla\opI[v_\delta](\gamma(t))\rangle_{L^2}\opdt\right| \leq \epsilon.\eqnum{\nexteqnno[ConstructingVPartI]}
$$
\endproclaim
\proclabel{ConstructingVPartI}
\proof Indeed, by Lemma \procref{BoundsAwayFromZero}, $\|\nabla\opI[v_\delta](\gamma(t))\|$ is uniformly bounded over $\Bbb{R}\setminus[-R,R]$. The result now follows since $g(t)$ decays exponentially as $t$ tends to $\pm\infty$.\qed
\proclaim{Lemma \nextprocno}
\noindent For all $R,\eta>0$,
$$
\mlim_{\delta\rightarrow 0}\int_{[-R,R]\setminus]-\eta,\eta[}\langle g(t),\tau(t)^{-1}\nabla\opI[v_\delta](\gamma(t))\rangle_{L^2}\opdt = 0.\eqnum{\nexteqnno[ConstructingVPartII]}
$$
\endproclaim
\proclabel{ConstructingVPartII}
\proof Since $u\in\Cal{U}_4$, the set of all $t\in[-R,R]\setminus]-\eta,\eta[$ such that $\gamma(t)$ contains $x_0$ is discrete and therefore finite. Its complement $\Omega$ thus has full measure in this set. However, for all $t\in\Omega$,
$$
\mlim_{\delta\rightarrow0}\|\nabla\opI[v_\delta](\gamma(t))\| = 0.
$$
The result now follows by Lemma \procref{BoundsAwayFromZero} and the dominated convergence theorem.\qed
\medskip
\noindent It thus remains to prove
\proclaim{Lemma \nextprocno}
\noindent There exists $c>0$ such that, for sufficiently small $\eta$,
$$
\mliminf_{\delta\rightarrow 0}\int_{-\eta}^\eta\langle g(t),\tau(t)^{-1}\nabla\opI[v_\delta](\gamma(t))\rangle_{L^2}\opdt\geq c.\eqnum{\nexteqnno[ConstructingVPartIII]}
$$
\endproclaim
\proclabel{ConstructingVPartIII}
\noindent Lemma \procref{ConstructingVPartIII} will follow from Lemma \procref{EstimateOverEtaInterval}, below.
\newsubhead{Technical properties of $v$}[TechnicalPropertiesOfV]
We continue to use the notation of the preceding sections. Trivially
$$
g_t(\theta,\phi) = g_0 + \opO(\left|\theta\right|+\left|\phi\right|+\left|t\right|).\eqnum{\nexteqnno[MorseSmaleAsymptoticsOfG]}
$$
Since $\tilde{v}_\delta$ is Lipschitz, it is almost everywhere differentiable. Its derivative satisfies
$$\eqalign{
d\tilde{v}_\delta &= \bigg(\frac{1}{\delta^4}\Bbb{I}_{0\leq t\leq\delta^3} - \frac{1}{\delta^4}\Bbb{I}_{\delta^3\leq t\leq 2\delta^3}\bigg)\chi^1_\delta(\theta)\chi^1_\delta(\phi)dt\cr
&\qquad\qquad\qquad\qquad + \opO\bigg(\frac{1}{\delta^2}\bigg)\Bbb{I}_{0\leq t\leq 2\delta^3}d\theta+ \opO\bigg(\frac{1}{\delta^2}\bigg)\Bbb{I}_{0\leq t\leq2\delta^3}d\phi.\cr}\eqnum{\nexteqnno[AsymptoticFormulaForIntegrand]}
$$
In the Fermi parametrisation $\Phi$, $N_t$ satisfies
$$
N_t = (1+\opO(t^2))\partial_t + \opO(t)\partial_\theta + \opO(t)\partial_\phi.\eqnum{\nexteqnno[MorseSmaleNormalField]}
$$
Since $T_0$ is minimal, the area form $\opdArea_t$ of the graph of $f_t$ satisfies
$$
\opdA_t = \frac{1}{2}(1+\opO(t^2))d\theta d\phi.\eqnum{\nexteqnno[MorseSmaleAreaForm]}
$$
It follows that, for any $C>0$, and for all $t\in[-C\delta,C\delta]$,
$$\eqalign{
g_t d\tilde{v}_\delta(N_t)dA_t &= \frac{g_0}{2\delta^4}(1 + \opO(\delta))\chi^1_\delta(\theta)\chi^1_\delta(\phi)\big(\Bbb{I}_{0\leq f_t\leq\delta^3}-\Bbb{I}_{\delta^3\leq f_t\leq 2\delta^3}\big)d\theta d\phi\cr
&\qquad\qquad\qquad\qquad+\opO\bigg(\frac{t}{\delta^2}\bigg)\Bbb{I}_{0\leq f_t\leq 2\delta^3}d\theta d\phi.\cr}\eqnum{\nexteqnno[MorseSmaleAsymptoticsOfIntegrand]}
$$
Let $\tilde{\phi}:]-\delta,\delta[\rightarrow]-\delta,\delta[$ be such that the intersection of $\opZ^\opreg_A$ with the square $[\delta,\delta]^2$ coincides with the graph of $\tilde{\phi}$ over the $\theta$-axis. By construction,
$$
\tilde{\phi}(0) = \tilde{\phi}'(0) = 0.\eqnum{\nexteqnno[PropertiesOfTildePhi]}
$$
Let $a,b:[-\delta,\delta]\rightarrow\Bbb{R}$ be such that
$$
f_t(\theta,\phi) = a(\theta)\psi t + \frac{1}{2} b(\theta)t^2 + \opO(t\psi^2 + t^2\psi + t^3),\eqnum{\nexteqnno[MorseSmaleAsymptoticsOfF]}
$$
where $\psi := \phi - \tilde{\phi}(\theta)$. Without loss of generality, we may suppose that
$$\multiline{
a_0:=a(0)>0\ \text{and}\cr
b_0:=b(0)<0.\cr}\eqnum{\nexteqnno[MorseSmaleAAndB]}
$$
Lemma \procref{ConstructingVPartIII} now follows from
\proclaim{Lemma \nextprocno}
\noindent For sufficiently small $\eta$,
$$
\mliminf_{\delta\rightarrow 0}\int_{-\eta}^\eta\int_{\Bbb{T}}g_t d\tilde{v}_\delta(N_t)\opdA_t\opdt\geq\frac{9g_0}{8a_0}.\eqnum{\nexteqnno[EstimateOverEtaInterval]}
$$
\endproclaim
\proclabel{EstimateOverEtaInterval}
We will estimate \eqnref{EstimateOverEtaInterval} by studying the preimages under $f_t$ of the sets $[0,\delta^3]$ and $[\delta^3,2\delta^3]$. First, since $\delta<\eta$, for $(\theta,\phi,t)\in[-\delta,\delta]^2\times[-\eta,\eta]$,
$$
\frac{\partial f_t}{\partial\phi}(\theta,\phi) = (a_0+\opO(\eta))t.\eqnum{\nexteqnno[MorseSmaleDerivativeOfF]}
$$
In particular, for sufficiently small $\eta$ and for all non-zero $t\in[-\eta,\eta]$, $f_t$ is strictly monotone in $\phi$, increasing when $t>0$ and decreasing when $t<0$. We now restrict attention to $t>0$, since the case where $t<0$ then follows by reversing the signs of $t$ and $\phi$. Define $\tilde{\phi}_{1,\delta},\tilde{\phi}_{2,\delta},\tilde{\phi}_{3,\delta}:[-\delta,\delta]\times]0,\eta]\rightarrow]-\infty,\delta]$ by
$$\eqalign{
\tilde{\phi}_{1,\delta}(\theta,t) &:= \msup\left\{\phi\in[-\delta,\delta]\ | f_t(\theta,\phi)\leq0\right\},\cr
\tilde{\phi}_{2,\delta}(\theta,t) &:= \msup\left\{\phi\in[-\delta,\delta]\ | f_t(\theta,\phi)\leq\delta^3\right\}\ \text{and}\cr
\tilde{\phi}_{3,\delta}(\theta,t) &:= \msup\left\{\phi\in[-\delta,\delta]\ | f_t(\theta,\phi)\leq2\delta^3\right\}.\cr}\eqnum{\nexteqnno[DefinitionOfTildePhi]}
$$
Observe that, for all $(\delta,\theta,t)$,
$$
\tilde{\phi}_{1,\delta}(\theta,t) \leq \tilde{\phi}_{2,\delta}(\theta,t)\leq \tilde{\phi}_{3,\delta}(\theta,t)\leq\delta.
\eqnum{\nexteqnno[FirstRelationOfPhis]}
$$
\proclaim{Lemma \nextprocno}
\noindent For sufficiently small $\eta$, there exists $B_1>0$ such that, for all $\delta<\eta$ and for all $t\in]0,\eta]$,
$$
\tilde{\phi}_{1,\delta}(\theta,t)\geq\tilde{\phi}(\theta)\geq B_1\delta^2.\eqnum{\nexteqnno[MSEstimateI]}
$$
\endproclaim
\proclabel{MSEstimateI}
\proof By \eqnref{MorseSmaleAsymptoticsOfF},
$$
f_t(\theta,\tilde{\phi}(\theta)) = \frac{1}{2}b_0(1+\opO(\left|\theta\right|+\eta))t^2<0.
$$
It follows that $\tilde{\phi}_{1,\delta}(\theta,t)\geq\tilde{\phi}(\theta)$, as desired.\qed
\proclaim{Lemma \nextprocno}
\noindent For sufficiently small $\eta$, there exists $C_2>0$ such that, for all $\delta<\eta$, for all $\theta\in[-\delta,\delta]$ and for $C_2\delta<t<\eta$,
$$
\tilde{\phi}_{1,\delta}(\theta,t) = \delta.\eqnum{\nexteqnno[MSEstimateII]}
$$
\endproclaim
\proclabel{MSEstimateII}
\proof Indeed, by \eqnref{MorseSmaleAsymptoticsOfF},
$$
f_t(\theta,\delta) = a_0\delta t(1+\opO(\eta)) + \frac{1}{2}b_0t^2(1+\opO(\eta)).
$$
Since $a_0>0$ and $b_0<0$, there exists $B_1>0$ such that
$$\eqalign{
f_t(\theta,\delta)
&\leq a_0\delta t(1+B_1\eta) + \frac{1}{2}b_0t^2(1-B_1\eta)\cr
&=a\delta t(1+B_1\eta)\bigg(1 - \frac{\left|b_0\right|t(1-B_1\eta)}{2a_0\delta(1+B_1\eta)}\bigg).\cr}
$$
It follows that, for sufficiently small $\eta$, there exists $C_2>0$ such that, for all $t>C_2\delta$,
$$
f_t(\theta,\delta)\leq 0,
$$
so that
$$
\tilde{\phi}_{1,\delta}(\theta,t) = \delta.
$$
This completes the proof.\qed
\proclaim{Lemma \nextprocno}
\noindent For sufficiently small $\eta$, there exists $B_2>0$ and $C_3>0$ such that, for all $\delta<\eta$, for all $\theta\in[-\delta,\delta]$ and for all $0<t<(\delta^2/a_0)(1-C_3\delta)$,
$$\eqalign{
\tilde{\phi}_{1,\delta}(\theta,t) &\leq B_2\delta^2\ \text{and}\cr
\tilde{\phi}_{2,\delta}(\theta,t) &=\delta.\cr}\eqnum{\nexteqnno[MSEstimateIII]}
$$
\endproclaim
\proclabel{MSEstimateIII}
\proof Indeed, by \eqnref{MorseSmaleAsymptoticsOfF}, for $\theta\in[-\delta,\delta]$ and for $0<t<(\delta^2/a_0)(1-C_3\delta)$,
$$\eqalign{
f_t(\theta,\delta)
&=a_0t\delta(1+\opO(\delta)) + \opO(\delta^4)\cr
&\leq\delta^3 - C_3\delta^4 + \opO(\delta^4).\cr}
$$
Thus, for $C_3$ sufficiently large, for all such $\theta$ and $t$,
$$
f_t(\theta,\delta)\leq\delta^3.
$$
It follows that
$$
\tilde{\phi}_{2,\delta}(\theta,t) = \delta.
$$
By \eqnref{MorseSmaleAsymptoticsOfF} again, for all $\theta$ and $t$,
$$
f_t(\theta,\tilde{\phi}(\theta)) = \frac{1}{2}b_0(1+\opO(\delta))t^2<0.
$$
Furthermore, for all $\phi\in[-\delta,\delta]$,
$$
\frac{\partial f_t}{\partial\phi}(\theta,\phi) = a_0(1+\opO(\delta))t.
$$
There therefore exists $A>0$ such that, for all $\theta$ and $t$,
$$
f_t(\theta,\tilde{\phi}(\theta) + (b_0/2a_0)(1+A\delta)t)>0,
$$
so that
$$
\tilde{\phi}_{1,\delta}(\theta,t)\leq\tilde{\phi}(\theta) + \frac{b_0}{2a_0}(1+A\delta)t \leq B_2\delta^2,
$$
for some $B_2>0$. This completes the proof.\qed
\proclaim{Lemma \nextprocno}
\noindent For sufficiently small $\eta$, with $C_2$ and $C_3$ as in Lemmas \procref{MSEstimateII} and \procref{MSEstimateIII} respectively, there exists $B_3>0$ such that, for all $\delta<\eta$, for all $\theta\in[-\delta,\delta]$ and for all $(\delta^2/a_0)(1-C_3\delta)\leq t\leq C_2\delta$,
$$\eqalign{
\tilde{\phi}_{3,\delta}(\theta,t) - \tilde{\phi}_{2,\delta}(\theta,t) &\leq \frac{\delta^3}{a_0t}(1+B_3\delta)\ \text{and}\cr
\tilde{\phi}_{2,\delta}(\theta,t) - \tilde{\phi}_{1,\delta}(\theta,t) &\leq \frac{\delta^3}{a_0t}(1+B_3\delta),\cr}\eqnum{\nexteqnno[MSEstimateIVA]}
$$
and, if $\tilde{\phi}_{2,\delta}(\theta,t)<\delta$,
$$
\tilde{\phi}_{2,\delta}(\theta,t) - \tilde{\phi}_{1,\delta}(\theta,t) \geq \frac{\delta^3}{a_0t}(1-B_3\delta).\eqnum{\nexteqnno[MSEstimateIVB]}
$$
\endproclaim
\proclabel{MSEstimateIV}
\proof Indeed, by \eqnref{MorseSmaleAsymptoticsOfF}, for $0<t\leq C_2\delta$ and for all $(\theta,\phi)\in[-\delta,\delta]^2$,
$$
\frac{\partial f_t}{\partial\phi}(\theta,\phi) = a_0(1+\opO(\delta))t.
$$
It follows that
$$\eqalign{
\tilde{\phi}_{3,\delta}(\theta,t) - \tilde{\phi}_{2,\delta}(\theta,t) &\leq \frac{\delta^3}{a_0t}(1+B_3\delta)\ \text{and}\cr
\tilde{\phi}_{2,\delta}(\theta,t) - \tilde{\phi}_{1,\delta}(\theta,t) &\leq \frac{\delta^3}{a_0t}(1+B_3\delta),\cr}
$$
for some $B_3>0$. If $\tilde{\phi}_{2,\delta}(\theta,t)<\delta$, then
$$
f_t(\theta,\tilde{\phi}_{2,\delta}(\theta,t)) = \delta^3,
$$
so that, upon increasing $B_3$ if necessary,
$$
\tilde{\phi}_{2,\delta}(\theta,t) - \tilde{\phi}_{1,\delta}(\theta,t) \geq \frac{\delta^3}{a_0t}(1-B_3\delta).
$$
This completes the proof.\qed
\proclaim{Lemma \nextprocno}
\noindent For sufficiently small $\eta$, with $C_3$ as in Lemma \procref{MSEstimateIII},
$$
\mlim_{\delta\rightarrow 0}\int_0^{\frac{\delta^2}{a_0}(1-C_3\delta)}\int_{\Bbb{T}}g_t d\tilde{v}_\delta(N_t)\opdA_t\opdt = \frac{9g_0}{16a_0}.\eqnum{\nexteqnno[MSEstimateV]}
$$
\endproclaim
\proclabel{MSEstimateV}
\proof Indeed, by \eqnref{MorseSmaleAsymptoticsOfF},
$$\multiline{
\int_0^{\frac{\delta^2}{a_0}(1-C_3\delta)}\int_{\Bbb{T}}g_t d\tilde{v}_\delta(N_t)\opdA_t\opdt\cr
\qquad\qquad=\int_0^{\frac{\delta^2}{a_0}(1-C_3\delta)}\int_{-\delta}^\delta\int_{\tilde{\phi}_{1,\delta}(\theta,t)}^\delta
\frac{g_0}{2\delta^4}(1+\opO(\delta))\chi_\delta^1(\theta)\chi_\delta^1(\phi) + \opO(t\delta^{-2})d\phi d\theta\opdt.\cr}
$$
Since
$$
-B_1\delta^2 \leq \tilde{\phi}_{1,\delta}(\theta,t) \leq B_2\delta^2,
$$
with $B_1$ and $B_2$ as in Lemmas \procref{MSEstimateI} and \procref{MSEstimateIII}, it follows that
$$
\int_0^{\frac{\delta^2}{a_0}(1-C_3\delta)}\int_{\Bbb{T}}g_t d\tilde{v}_\delta(N_t)\opdA_t\opdt = \frac{9g_0}{16a_0}(1+\opO(\delta)).
$$
The result follows upon letting $\delta$ tend to zero.\qed
\proclaim{Lemma \nextprocno}
\noindent For sufficiently small $\eta$, with $C_2$ and $C_3$ as in Lemmas \procref{MSEstimateII} and \procref{MSEstimateIII} respectively,
$$
\mliminf_{\delta\rightarrow 0}\int_{\frac{\delta^2}{a_0}(1-C_3\delta)}^{C_2\delta}\int_{S_1\times S^1}g_t d\tilde{v}_\delta(N_t)\opdt \geq 0.\eqnum{\nexteqnno[MSEstimateVI]}
$$
\endproclaim
\proclabel{MSEstimateVI}
\proof Indeed, by \eqnref{MorseSmaleAsymptoticsOfF},
$$\multiline{
\int_{\frac{\delta^2}{a_0}(1-C_3\delta)}^{C_2\delta}\int_{\Bbb{T}}g_t d\tilde{v}_\delta(N_t)\opdA_t dt\cr
\qquad\qquad=\int_{\frac{\delta^2}{a_0}(1-C_3\delta)}^{C_2\delta}\int_{-\delta}^\delta\int_{\tilde{\phi}_{1,\delta}(\theta,t)}^{\tilde{\phi}_{2,\delta}(\theta,t)}\frac{g_0}{2\delta^4}(1+\opO(\delta))
\chi_\delta^1(\theta)\chi_\delta^1(\phi)d\phi d\theta\opdt\cr
\qquad\qquad\qquad\qquad-\int_{\frac{\delta^2}{a_0}(1-C_3\delta)}^{C_2\delta}\int_{-\delta}^\delta\int_{\tilde{\phi}_{2,\delta}(\theta,t)}^{\tilde{\phi}_{3,\delta}(\theta,t)}\frac{g_0}{2\delta^4}(1+\opO(\delta))
\chi_\delta^1(\theta)\chi_\delta^1(\phi)d\phi d\theta\opdt\cr
\qquad\qquad\qquad\qquad\qquad\qquad+\int_{\frac{\delta^2}{a_0}(1-C_3\delta)}^{C_2\delta}\int_{-\delta}^\delta\int_{\tilde{\phi}_{1,\delta}(\theta,t)}^{\tilde{\phi}_{3,\delta}(\theta,t)}\opO(t\delta^{-2})d\phi d\theta\opdt.\cr}
$$
Since $\tilde{\phi}_{1,\delta}(\theta,t)\geq\opO(\delta^2)$, $\chi_\delta^1(\phi)$ is monotone non-increasing as $\phi$ varies in the interval $[\tilde{\phi}_{1,\delta}(\theta,t),\tilde{\phi}_{3,\delta}(\theta,t)]$. It follows by Lemma \procref{MSEstimateIV} that, for all $(\delta^2/a_0)(1-C_3\delta)\leq t\leq C_2\delta$ and for all $\theta\in[-\delta,\delta]$,
$$
\int_{\tilde{\phi}_{1,\delta}(\theta,t)}^{\tilde{\phi}_{2,\delta}(\theta,t)}\frac{g_0}{2\delta^4}(1+\opO(\delta))\chi_\delta^1(\phi)d\phi
-\int_{\tilde{\phi}_{2,\delta}(\theta,t)}^{\tilde{\phi}_{3,\delta}(\theta,t)}\frac{g_0}{2\delta^4}(1+\opO(\delta))\chi_\delta^1(\phi)d\phi
\geq\frac{\opO(1)}{t}.
$$
Likewise, for all such $t$ and $\theta$,
$$
\int_{\tilde{\phi}_{1,\delta}(\theta,t)}^{\tilde{\phi}_{3,\delta}(\theta,t)}\opO(t\delta^{-2})d\phi\geq\opO(\delta).
$$
Combining these relations yields
$$
\int_{\frac{\delta^2}{a_0}(1-C_3\delta)}^{C_2\delta}\int_{\Bbb{T}}g_t d\tilde{v}_\delta(N_t)\opdA_t dt
\geq \int_{\frac{\delta^2}{a_0}(1-C_3\delta)}^{C_2\delta}\int_{-\delta}^\delta\frac{\opO(1)}{t}d\theta\opdt
=\opO(\delta\opLog(\delta)).
$$
The result follows upon letting $\delta$ tend to zero. This completes the proof.\qed
\medskip
\noindent We now prove Lemma \procref{EstimateOverEtaInterval}.
\medskip
{\bf\noindent Proof of Lemma \procref{EstimateOverEtaInterval}:\ } By Lemmas \procref{MSEstimateV} and \procref{MSEstimateVI},
$$
\mliminf_{\delta\rightarrow 0}\int_0^{C_1\delta}\int_{\Bbb{T}}g_t d\tilde{v}_\delta(N_t)\opdA_t\opdt \geq \frac{9g_0}{16a_0}.
$$
Upon reversing the signs of $t$ and $\phi$, these lemmas also yield
$$
\mliminf_{\delta\rightarrow 0}\int_{-C_1\delta}^0\int_{\Bbb{T}}g_t d\tilde{v}_\delta(N_t)\opdA_t\opdt \geq \frac{9g_0}{16a_0}.
$$
Finally, by Lemma \procref{MSEstimateII},
$$
\int_{[-\eta,\eta]\setminus[-C_1\delta,C_1\delta]}\int_{\Bbb{T}}g_t d\tilde{v}_\delta(N_t)\opdA_t\opdt = 0.
$$
The result now follows upon combining these relations.\qed
\newhead{Perturbation theory}[PerturbationTheory]
\newsubhead{White's Construction}[WhitesConstruction]
We now reformulate the construction described by White in Section $3$ of \cite{White} in terms of our current framework. For all $(k,\alpha)$ with $k\geq 1$, let $\Cal{G}^{1,\alpha}(\Bbb{S}^3)$ be the Banach manifold of $C^{k,\alpha}$ riemannian metrics over $\Bbb{S}^3$. Define
$$
\opH:\opO(4)\times C^{2,\alpha}(\Bbb{T},]-\pi/4,\pi/4[)\times\Cal{G}^{1,\alpha}(\Bbb{S}^3)\rightarrow C^{0,\alpha}(\Bbb{T})
$$
such that, for all $M\in\opO(4)$, for all $f\in C^{2,\alpha}(\Bbb{T},]-\pi/4,\pi/4[)$, for all $g\in\Cal{G}^{1,\alpha}(\Bbb{S}^3)$ and for all $(\theta,\phi)\in\Bbb{T}$, $\opH[M,f,g](\theta,\phi)$ is the mean curvature of the embedding $M\ope[f]$ with respect to the metric $g$ at the point $(\theta,\phi)$. Observe that $\opH$ is a smooth function between Banach manifolds.
\par
Let $g_1$ be the standard metric of $\Bbb{S}^3$. For all $M\in\opO(4)$, the partial derivative of $\opH$ with respect to the second component at the point $(M,0,g_1)$ is given by
$$
\opD_2\opH[M,0,g_1]f = -2(\Delta+2)f.\eqnum{\nexteqnno[DerivativeOfMeanCurvatureFunctional]}
$$
Recall that the kernel $\opK$ of $(\Delta+2)$ in $L^2(\Bbb{T})$ consists of all functions of the form $\psi_A$ for some $A\in\opEnd(\Bbb{R}^2)$. The fact that this kernel is non-trivial means that Clifford tori do not, in general, perturb to minimal surfaces for metrics close to $g_1$. White's construction compensates for this degenericity by instead perturbing Clifford tori to surfaces whose mean curvature is as small as possible in a certain algebraic sense which we will make clear presently. We proceed as follows. Let $\opK^\perp$ be the orthogonal complement of $K$ in $L^2(\Bbb{T})$. Let $\pi:L^2(\Bbb{T})\rightarrow\opK$ and $\pi^\perp:L^2(\Bbb{T})\rightarrow\opK^\perp$ be the orthogonal projections. For all $(k,\alpha)$, let $\opK^{\perp,k,\alpha}$ be the intersection of $\opK^\perp$ with $C^{k,\alpha}(\Bbb{T})$. By classical PDE theory (see \cite{Taylor}), the operator
$$
\pi^\perp\circ(\Delta+2)
$$
defines a continuous linear isomorphism from $K^{\perp,2,\alpha}$ into $K^{\perp,0,\alpha}$. By the inverse function theorem, there therefore exists a neighbourhood $\Omega^{k,\alpha}$ of $g_1$ in $\Cal{G}^{k,\alpha}(\Bbb{S}^3)$ having the property that there exists a unique smooth function $\optildef:\opO(4)\times\Omega^{k,\alpha}\rightarrow\opK^{\perp,2,\alpha}$ such that, for all $M\in\opO(4)$,
$$
\optildef[M,g_1] = 0\eqnum{\nexteqnno[DefinitionOfFI]}
$$
and, for all $g\in\Omega^{k,\alpha}$,
$$
(\pi^\perp\circ\opH)[M,\optildef[M,g],g] = 0.\eqnum{\nexteqnno[DefinitionOfFII]}
$$
In other words, the function $\optildef$ has been defined in such a manner that, for all $(M,g)\in\opO(4)\times\Omega^{k,\alpha}$, the mean curvature of the embedding $M\ope[\optildef[M,g]]$ with respect to $g$ vanishes modulo the finite-dimensional subspace $K$. We henceforth refer to $\optildef$ as {\sl White's functional}.
\par
Define $\optildeh:\opO(4)\times\Omega^{k,\alpha}\rightarrow K$ by
$$
\optildeh[M,g] := (\pi\circ\opH)[M,\optildef[M,g],g],\eqnum{\nexteqnno[DefinitionOfH]}
$$
so that, for all $(M,g)\in\opO(4)\times\Omega^{k,\alpha}$ and for all $(\theta,\phi)\in\Bbb{T}$, $\optildeh[M,g](\theta,\phi)$ is the mean curvature of the embedding $M\ope[\optildef[M,g]]$ with respect to the metric $g$ at the point $(\theta,\phi)$. Define $\optildea:\opO(4)\times\Omega^{k,\alpha}\rightarrow\Bbb{R}$ such that, for all $M\in\opO(4)$ and for all $g\in\Omega$, $\optildea[M,g]$ is the area of $M\ope[\optildef[M,g]]$ with respect to the metric $g$. In what follows, we will only be concerned with smooth, finite-dimensional families of smooth metrics about $g_1$. We therefore suppress the superscript $(k,\alpha)$ and write $\Cal{G}(\Bbb{S}^3)$ instead of $\Cal{G}^{k,\alpha}(\Bbb{S}^3)$ and $\Omega$ instead of $\Omega^{k,\alpha}$.
\par
By uniqueness, the functions $\optildef$, $\optildeh$ and $\optildea$ are invariant under certain transformations. First, since $\opO(4)$ is compact, we may suppose that $\Omega$ is invariant under the pull-back action of this group on $\Cal{G}(\Bbb{S}^3)$. Now, for all $M,N\in\opO(4)$, for all $f\in C^{2,\alpha}(\Bbb{T},]-\pi/4,\pi/4[)$ and for all $g\in\Omega$,
$$
\optildeh[M,f,N^*g] = \optildeh[NM,f,g],
$$
so that, by uniqueness, for all $M,N\in\opO(4)$ and for all $g\in\Omega$,
$$\eqalign{
\opf[M,N^*g] &= \opf[NM,g],\cr
\optildeh[M,N^*g] &= \optildeh[NM,g]\ \text{and}\cr
\optildea[M,N^*g] &= \optildea[NM,g].\cr}\eqnum{\nexteqnno[LeftInvariance]}
$$
Likewise, for all $M\in\opO(4)$, for all $N\in\opO(2)^2$, for all $f\in C^{2,\alpha}(\Bbb{T},]-\pi/4,\pi/4[)$ and for all $g\in\Omega$,
$$
\optildeh[MN,f\circ N,g] = \optildeh[M,f,g]\circ N,
$$
so that, by uniqueness again, for all $M\in\opO(4)$, for all $N\in\opO(2)^2$ and for all $g\in\Omega$,
$$\eqalign{
\opf[MN,g] &= \opf[M,g]\circ N,\cr
\optildeh[MN,g] &= \optildeh[M,g]\circ N\ \text{and}\cr
\optildea[MN,g] &= \optildea[M,g].\cr}\eqnum{\nexteqnno[RightInvariance]}
$$
\newsubhead{The derivatives of $\opa$}[TheDerivativesOfA]
We now determine the derivatives of the function $\opa$ defined in the preceding section. For each $i$, let $\opD_i\opa$ be its partial derivative with respect to the $i$'th component.
\proclaim{Lemma \nextprocno}
\noindent If $h=2ug_1$ for some $u\in C^\infty(\Bbb{S}^3)$ then, for all $M\in\opO(4)$,
$$
\opD_2\optildea[M,g_1]h = (\opI[u]\circ\opC)(M),\eqnum{\nexteqnno[FirstDerivativeOfTildeAWRTG]}
$$
where $\opI$ and $\opC$ are the functions defined in \eqnref{DefinitionOfFunctionalI} and \eqnref{DiffeoOfOWithCL} respectively.
\endproclaim
\proclabel{FirstDerivativeOfTildeAII}
\proof By the left invariance of \eqnref{LeftInvariance}, we may suppose that $M=\opId$. Let $(g_{1+s})_{s\in]-\epsilon,+\epsilon[}$ be a smooth family of metrics about $g_1$ such that
$$
\frac{\partial g_{1+s}}{\partial s}\bigg|_{s=0} = h = 2ug_1.
$$
For sufficiently small $s$, denote $e_s(\theta,\phi):=(\theta,\phi,\optildef[\opId,g_s](\theta,\phi))$. For sufficiently small $s$ and $t$, let $a_{s,t}$ be the area of the embedding $e_s$ with respect to the metric $\Phi^*g_{1+t}$. Since $T_0$ is minimal, the area form of $e_s$ with respect to $\Phi^*g_1$ satisfies
$$
\opdA_{s,0} := \frac{1}{2}d\theta d\phi + \opO(s^2),
$$
so that
$$
\frac{\partial a_{s,t}}{\partial s}\bigg|_{s,t=0}
= \frac{1}{2}\frac{\partial}{\partial s}\int_{\Bbb{T}}d\theta d\phi\bigg|_{s=0}=0.
$$
On the other hand, the area form of $e_0$ with respect to $\Phi^*g_{1+t}$ is
$$
\opdA_{0,t} := \frac{1}{2}(1+2ut)d\theta d\phi + \opO(t^2),
$$
so that
$$
\frac{\partial a_{s,t}}{\partial t}\bigg|_{s,t=0}
= \frac{1}{2}\frac{\partial}{\partial t}\int_{\Bbb{T}}(1+2ut)d\theta d\phi\bigg|_{t=0} = \opI[u].
$$
Applying the chain rule now yields
$$
\opD_2\optildea[\opId,g_1]h = \frac{\partial a_{s,s}}{\partial s}\bigg|_{s=0} = \frac{\partial a_{s,t}}{\partial s}\bigg|_{s,t=0} + \frac{\partial a_{s,t}}{\partial t}\bigg|_{s,t=0} = \opI[u](T_0),
$$
as desired.\qed
\proclaim{Lemma \nextprocno}
\noindent For all $M\in\opO(4)$, for all $g\in\Omega$ and for all $A\in\opEnd(\Bbb{R}^2)$,
$$
\opD_1\optildea[M,g](\tau^L\xi_A)(M) = \frac{1}{2}\int_{\Bbb{T}}\psi_A\optildeh[M,g]d\theta d\phi + \opO\big(\|A\|\|g-g_1\|_{C^{1,\alpha}}^2\big),\eqnum{\nexteqnno[FirstDerivativeOfTildeA]}
$$
where $\tau^L\xi_A$ here denotes the left-translation of $\xi_A$ over $\opO(4)$.

\endproclaim
\proclabel{FirstDerivativeOfTildeA}
\remark By the right-invariance of \eqnref{RightInvariance}, $\opD_1\optildea[M,g]$ vanishes over $(\tau^L\frak{h})(M)$ so that \eqnref{FirstDerivativeOfTildeA} completely describes $\opD_1\optildea[M,g]$.
\medskip
\proof By the left invariance of \eqnref{LeftInvariance}, we may suppose that $M=\opId$. Observe now that, since $\optildeh[\opId,g_1]=0$,
$$
\optildeh[\opId,g] = \opO\big(\|g-g_1\|_{C^{1,\alpha}}\big).\eqnum{\nexteqnno[TildeHIsSmall]}
$$
Now choose $A\in\opEnd(\Bbb{R}^2)$. For all sufficiently small $s$ and $t$, let $a_{s,t}$ be the area of the embedding $\opExp(s\xi_A)\ope[\optildef[\opExp(t\xi_A),g]]$ with respect to the metric $g$. Let $X_A$ and $F_A$ be as in Section \subheadref{KillingFieldsOverCL}. For all sufficiently small $s$ and $t$, denote
$$
e_{s,t}(\theta,\phi) := F_{A,s}(\theta,\phi,\optildef[\opExp(t\xi_A),g]).
$$
Then
$$\eqalign{
(\Phi\circ e_{s,t})(\theta,\phi) &= (\Phi\circ F_{A,s})(\theta,\phi,\optildef[\opExp(t\xi_A),g])\cr
&= \opExp(s\xi_A)\ope[\optildef[\opExp(t\xi_A),g]](\theta,\phi),\cr}
$$
so that $a_{s,t}$ is equal to the area of $e_{s,t}$ with respect to the metric $\Phi^*g$. By definition of $X_A$ and $F_A$,
$$
\frac{\partial e_{s,t}}{\partial s}(\theta,\phi)\bigg|_{s,t=0} = X_A(\theta,\phi,\optildef[\opId,g]).
$$
Thus, since the unit normal vector  field over $e_{0,0}$ with respect to $\Phi^*g$ is
$$
N := (0,0,1) + \opO\big(\|g-g_1\|_{C^{1,\alpha}}\big),
$$
and since
$$
\optildef[\opId,g] = \opO(\|g-g_1\|_{C^{1,\alpha}}),
$$
it follows by \eqnref{ExplicitFormulaForX} that
$$
(\Phi^*g)\bigg(\frac{\partial e_{s,t}}{\partial s}(\theta,\phi)\bigg|_{s,t=0},N(\theta,\phi)\bigg) = \psi_A(\theta,\phi) + \opO\big(\|A\|\|g-g_1\|_{C^{1,\alpha}}\big).
$$
Since the area form of $e_{0,0}$ with respect to $\Phi^*g$ is
$$
\opdA = \frac{1}{2}d\theta d\phi + \opO\big(\|g-g_1\|_{C^{1,\alpha}}),
$$
it follows by \eqnref{TildeHIsSmall} and the first variation formula for area (c.f. Equation $1.45$ of \cite{ColdMinII}) that
$$
\frac{\partial a_{s,t}}{\partial s}\bigg|_{s,t=0} = \frac{1}{2}\int_{\Bbb{T}}\optildeh[\opId,g]\psi_A(\theta,\phi)+\opO\big(\|A\|\|g-g_1\|^2_{C^{1,\alpha}}\big)d\theta d\phi.
$$
On the other hand, for all $t$,
$$
\optildef[\opExp(t\xi_A),g_1] = 0,
$$
so that
$$
\frac{\partial}{\partial t}\optildef[\opExp(t\xi_A),g] = \opO\big(\|A\|\|g-g_1\|_{C^{1,\alpha}}\big).
$$
It follows by \eqnref{TildeHIsSmall} and the first variation formula for area again that
$$
\frac{\partial a_{s,t}}{\partial t}\bigg|_{s,t=0} = \opO\big(\|A\|\|g-g_1\|^2_{C^{1,\alpha}}\big).
$$
Applying the chain rule now yields
$$
D_1\opa[T_0,g]\xi_A
= \frac{1}{2}\int_{\Bbb{T}}\optildeh[\opId,g]\psi_A(\theta,\phi)d\theta d\phi+\opO\big(\|A\|\|g-g_1\|^2_{C^{1,\alpha}}\big),
$$
as desired.\qed
\proclaim{Lemma \nextprocno}
\noindent Let $\opHess_1(\optildea)$ be the Hessian of $\optildea$ with respect to the first component. For all $M\in\opO(4)$, for all $g\in\Omega$ and for all $A,B\in\opEnd(\Bbb{R}^2)$,
$$\multiline{
\opHess_1(\optildea)[M,g]((\tau^L\xi_A)(M),(\tau^L\xi_B)(M))\cr
\qquad\qquad\qquad = \frac{1}{2}\int_{\Bbb{T}}\psi_A\opD_1\optildeh[M,g](\tau^L\xi_B)(M)d\theta d\phi + \opO\big(\|A\|\|B\|\|g-g_1\|^2_{C^{1,\alpha}}\big),\cr}\eqnum{\nexteqnno[HessianOfArea]}
$$
where $\tau^L\xi_A$ and $\tau^L\xi_B$ here denote respectively the left-translations of $\xi_A$ and $\xi_B$ over $\opO(4)$.
\endproclaim
\proclabel{HessianOfArea}
\proof By the left invariance of \eqnref{LeftInvariance}, we may suppose that $M=\opId$. Choose $A,B\in\opEnd(\Bbb{R}^2)$. Recall that
$$
\opExp(s\xi_A + t\xi_B) = \opExp(t\xi_B)\opExp(s\xi_A)\opExp((1/2)st[\xi_A,\xi_B]) + \opO(s^2 + t^2),
$$
so that, since $[\xi_A,\xi_B]\in\frak{h}$,
$$
\opExp(s\xi_A + t\xi_B)T_0 = \opExp(t\xi_B)\opExp(s\xi_A)T_0 + \opO(s^2 + t^2).
$$
Thus, bearing in mind Lemma \procref{FirstDerivativeOfTildeA},
$$\eqalign{
\opHess_1(\optildea)[\opId,g](\xi_A,\xi_B)
&=\frac{\partial^2}{\partial s\partial t}\optildea[\opExp(s\xi_A+t\xi_B),g]\bigg|_{s,t=0}\cr
&=\frac{\partial^2}{\partial s\partial t}\optildea[\opExp(t\xi_B)\opExp(s\xi_A),g]\bigg|_{s,t=0}\cr
&=\frac{\partial^2}{\partial s\partial t}\optildea[\opExp(s\xi_A),\opExp(t\xi_B)^*g]\bigg|_{s,t=0}\cr
&=\frac{\partial}{\partial t}\frac{1}{2}\int_{\Bbb{T}}\psi_A\optildeh[\opId,\opExp(t\xi_B)^*g]d\theta d\phi\bigg|_{t=0} + \opO(\|A\|\|B\|\|g-g_1\|^2_{C^{1,\alpha}})\cr
&=\frac{1}{2}\int_{\Bbb{T}}\psi_A\opD_1\optildeh[\opId,g]\xi_B d\theta d\phi + \opO(\|A\|\|B\|\|g-g_1\|^2_{C^{1,\alpha}}),\cr}
$$
as desired.\qed
\newsubhead{The mean curvature flow operator}[TheMeanCurvatureFlowOperator]
We now prove the main result of this paper. It remains only to express the existence of eternal mean curvature flows in terms of the vanishing of some functional, and to prove that zeroes of this functional can be constructed perturbatively. We will use the formalism of anisotropic H\"older spaces, which for the readers convenience we review in Appendix \headref{HoelderSpaces}.
\par
Fix $u\in C^\infty(\Bbb{S}^3)$, let $g_1$ be the standard metric of $\Bbb{S}^3$ and, for sufficiently small $\epsilon$, define the smooth family of metrics $(g_{1+s})_{s\in]-\epsilon,+\epsilon[}$ by
$$
g_{1+s} := e^{2su}g_1.\eqnum{\nexteqnno[DerivativeOfFamilyOfMetrics]}
$$
Let $\gamma:\Bbb{R}\rightarrow\opCL$ be a complete gradient flow of $\opI[u]$. Upon applying an element of $\opO(4)$. We may suppose that $\gamma(0)=T_0$. Let $M:\Bbb{R}\rightarrow\opO(4)$ be the unique lift of $\gamma$ such that $M(0)=\opId$ and such that, for all $t$,
$$
M(t)^{-1}\dot{M}(t) \in \frak{k}.\eqnum{\nexteqnno[DefinitioOfLiftOfGamma]}
$$
For $\eta\in C^{1,\alpha}(\Bbb{R},\frak{k})$, for $f\in C^{1,\alpha}_\opin(\Bbb{T}\times\Bbb{R},]-\pi/4,\pi/4[)$ and for $s\in]-\epsilon,\epsilon[$, define the flow $\ope[\eta,f,s]$ such that, for all $(\theta,\phi,t)\in \Bbb{T}\times\Bbb{R}$,
$$
\ope[\eta,f,s](\theta,\phi,t):=\opExp(\eta(st))M(st)\Phi(\theta,\phi,\optildef[\opExp(\eta(st))M(st),g_{1+s}](\theta,\phi)+sf(\theta,\phi,t)),
$$
where $\optildef$ is White's functional, as defined in Section \subheadref{WhitesConstruction}. For all such $\eta$, $f$ and $s$, define $\opN[\eta,f,s]$ and $\opH[\eta,f,s]$ such that, for all $(\theta,\phi,t)\in\Bbb{T}\times\Bbb{R}$, $\opN[\eta,f,s](\theta,\phi,t)$ and $\opH[\eta,f,s](\theta,\phi,t)$ are respectively the unit normal vector and mean curvature of the embedding $\ope[\eta,f,s](\cdot,\cdot,t)$ with respect to the metric $g_{1+s}$ at the point $(\theta,\phi)$. The {\sl mean curvature flow operator} is now defined by
$$
\opMCF[\eta,f,s] := g_{1+s}\bigg(\frac{\partial}{\partial t}\ope[\eta,f,s],\opN[\eta,f,s]\bigg) + \opH[\eta,f,s].\eqnum{\nexteqnno[DefinitionOfMCF]}
$$
This operator defines a smooth function from the Banach manifold $C^{1,\alpha}(\Bbb{R},\frak{k})\times C^{1,\alpha}_\opin(\Bbb{T}\times\Bbb{R},]-\pi/4,\pi/4[)\times]-\epsilon,\epsilon[$ into the Banach space $C_\opin^{0,\alpha}(\Bbb{T}\times\Bbb{R})$ which vanishes if and only if $\ope[\eta,f,s]$ is an eternal mean curvature flow with respect to the metric $g_{1+s}$.
\par
We now determine an explicit formula for $\opMCF$. First, define $\Psi:\frak{k}\rightarrow K$ such that, for all $A\in\opEnd(\Bbb{R}^2)$,
$$
\Psi(\xi_A) = \psi_A.
$$
\proclaim{Lemma \nextprocno}
\noindent For all $(\eta,f,s)$ and for all $(\theta,\phi,t)$,
$$\multiline{
g_{1+s}\bigg(\frac{\partial}{\partial t}\ope[\eta,f,s](\theta,\phi,t),\opN[\eta,f,s](\theta,\phi,t)\bigg)\cr
\qquad\qquad=s\Psi(\dot{\eta}(t) + M^{-1}(t)\dot{M}(t))(\theta,\phi) + s\frac{\partial f}{\partial t}(\theta,\phi,t) + \opO(s^2 + s\|\eta\|_{C^{1,\alpha}}^2).\cr}\eqnum{\nexteqnno[AsymptoticFormulaForNormalVariation]}
$$
\endproclaim
\proclabel{FirstComponentOfMCFOp}
\proof By time invariance of the mean curvature flow operator, we may suppose that $t=0$. Let $A,B\in\opEnd(\Bbb{R}^2)$ be such that
$$\eqalign{
\dot{\eta}(0) &= \xi_A\ \text{and}\cr
\dot{M}(0) &= \xi_B.\cr}
$$
Let $X_A$, $X_B$, $F_A$ and $F_B$ be as in Section \subheadref{KillingFieldsOverCL}. Define $\tilde{e}$ by
$$
\tilde{e}(\theta,\phi,t) = (F_{A,st}\circ F_{B,st})(\theta,\phi,\optildef[\opExp(\eta(st))M(st),g_{1+s}](\theta,\phi) + sf(\theta,\phi,t)).
$$
Then
$$\eqalign{
(\Phi\circ\tilde{e})(\theta,\phi,t) &= \opExp(st\xi_A)\opExp(st\xi_B)\Phi(\theta,\phi,\optildef[\opExp(\eta(st))M(st),g_{1+s}](\theta,\phi)+sf(\theta,\phi,t))\cr
&=\opExp(\eta(0))^{-1}\ope[\eta,f,s](\theta,\phi,t) + \opO(s^2t^2 + st\|\eta\|_{C^{1,\alpha}}^2).\cr}
$$
However,
$$\eqalign{
\frac{\partial\tilde{e}}{\partial t}(\theta,\phi,t)\bigg|_{t=0}
&=sX_A(\theta,\phi,\optildef[\opExp(\eta(0)),g_{1+s}](\theta,\phi) + sf(\theta,\phi,0))\cr
&\qquad + sX_B(\theta,\phi,\optildef[\opExp(\eta(0)),g_{1+s}](\theta,\phi) + sf(\theta,\phi,0))\cr
&\qquad\qquad + \bigg(0,0,s\frac{\partial}{\partial t}\optildef[\opExp(\eta(t))M(t),g_{1+s}](\theta,\phi)\bigg|_{t=0}\bigg)\cr
&\qquad\qquad\qquad + \bigg(0,0,s\frac{\partial f}{\partial t}(\theta,\phi,0)\bigg)+\opO(s\|\eta\|_{C^{1,\alpha}}^2).\cr}
$$
Since
$$
\optildef[\opExp(\eta(t))M(t),g_{1+s}] = \opO(s).
$$
It follows by \eqnref{ExplicitFormulaForX} that
$$
\frac{\partial\tilde{e}}{\partial t}(\theta,\phi,t)\bigg|_{t=0}
=\bigg(0,0,s\psi_A(\theta,\phi) + s\psi_B(\theta,\phi) + s\frac{\partial f}{\partial t}(\theta,\phi,0)\bigg) + \opO(s^2 + s\|\eta\|_{C^{1,\alpha}}^2).
$$
Since the unit normal vector field over $\tilde{e}(\cdot,\cdot,0)$ with respect to $g_{1+s}$ satisfies
$$
\tilde{N}(\theta,\phi) = (0,0,1) + \opO(s),
$$
it follows that
$$\multiline{
g_{1+s}\bigg(\frac{\partial}{\partial t}\ope[\eta,f,s](\theta,\phi,t),\optildeN[\eta,f,s](\theta,\phi,t)\bigg)\cr
\qquad\qquad=\big(\opExp(\eta(0))\circ\Phi)^*g_{1+s}\bigg(\frac{\partial\tilde{e}}{\partial t}(\theta,\phi,t)\bigg|_{t=0},\tilde{N}(\theta,\phi)\bigg)\cr
\qquad\qquad=s\psi_A(\theta,\phi)+s\psi_B(\theta,\phi) + s\frac{\partial f}{\partial t}(\theta,\phi,0) + \opO(s^2 + s\|\eta\|_{C^{1,\alpha}}^2)\cr
\qquad\qquad=s\Psi(\dot{\eta}(0) + M^{-1}(0)\dot{M}(0))(\theta,\phi) + s\frac{\partial f}{\partial t}(\theta,\phi) + \opO(s^2 + s\|\eta\|_{C^{1,\alpha}}^2),\cr}
$$
as desired.\qed
\medskip
\noindent Before determining the formula for $\opH[\eta,f,s]$, we first clarify some notation. Let $\opI$ and $\opC$ be the functions given by \eqnref{DefinitionOfFunctionalI} and \eqnref{DiffeoOfOWithCL} respectively, and let $\Psi$ be as above. For all $M\in\opO(4)$, $\nabla(\opI[u]\circ\opC)(M)$ is a tangent vector to $\opO(4)$ at the point $M$. We denote by $M^{-1}\nabla(\opI[u]\circ\opC)(M)$ the left translation of this vector to $\frak{o}(4)$. By the invariance of \eqnref{LeftInvariance}, this left-translated vector is an element of $\frak{k}$, and thus lies in the domain of $\Psi$. Likewise, for all $M\in\opO(4)$ and for all $\xi\in\frak{k}$, the vector $M^{-1}\opHess(\opI[u]\circ\opC)(M)\tau^L\xi(M)$ also lies in the domain of $\Psi$. We obtain
\proclaim{Lemma \nextprocno}
\noindent For all $\eta$, for all $f$ and for all $s$,
$$\eqalign{
\opH[\eta,f,s](\theta,\phi,t)
&=s\Psi M(t)^{-1}\nabla(\opI[u]\circ\opC)(M(t))\cr
&\qquad\qquad + s\Psi M(t)^{-1}\opHess(\opI[u]\circ\opC)(M(t))\tau^L\eta(t)(M(t))\cr
&\qquad\qquad\qquad\qquad - 2s(\Delta+2)f(\theta,\phi,t) + \opO(s^2 + s\|\eta\|_{C^{1,\alpha}}^2),\cr}\eqnum{\nexteqnno[AsymptoticFormulaForMeanCurvature]}
$$
where $\opC$ is the function defined in \eqnref{DiffeoOfOWithCL}.

\endproclaim
\proclabel{SecondComponentOfMCFOp}
\proof By the invariance of \eqnref{LeftInvariance}, we may suppose that $\gamma(0)=T_0$ and $M(0)=\opId$. Likewise, by the time-invariance of the mean curvature flow operator, we may suppose that $t=0$. Denote
$$
e_0[\eta,f,s](\theta,\phi) := \opExp(\eta(0))\Phi(\theta,\phi,\opf[\opExp(\eta(0)),g_{1+s}]).
$$
By construction, the mean curvature of $e_0$ with respect to the metric $g_{1+s}$ is equal to $\optildeh[\opExp(\eta(0)),g_{1+s}]$. The Jacobi operator of $\ope_0[\eta,f,0]$ with respect to the metric $g_1$ is
$$
J_0 := -2(\Delta+2),
$$
so that the Jacobi operator of $e_0[\eta,f,0]$ with respect to the metric $g_{1,s}$ is
$$
J_s := -2(\Delta+2) + \opO(s).
$$
The mean curvature of $\ope[\eta,f,s]$ is thus
$$
\opH[\eta,f,s](\theta,\phi,0)
=\optildeh[\opExp(\eta(0)),g_{1+s}](\theta,\phi) - 2s(\Delta+2)f(\theta,\phi,0) + \opO(s^2).
$$
The result now follows by Lemmas \procref{FirstDerivativeOfTildeAII}, \procref{FirstDerivativeOfTildeA} and \procref{HessianOfArea}.\qed
\medskip
\noindent Combining these relations yields
\proclaim{Lemma \nextprocno}
\noindent For all $\eta$, for all $f$ and for all $s$
$$\eqalign{
\opMCF[\eta,f,s] &= s\Psi\bigg(\frac{\partial\eta}{\partial t}(st) + M(st)^{-1}\opHess(\opI[u]\circ\opC)(M(st))\tau^L\eta(st)(M(st))\bigg)\cr
&\qquad\qquad + s\bigg(\frac{\partial f}{\partial t} - 2(\Delta+2)f\bigg) + \opO(s^2 + s\|\eta\|_{C^{1,\alpha}}^2).\cr}\eqnum{\nexteqnno[AsymptoticFormulaForMeanCurvatureFlowOperator]}
$$
\endproclaim
\proclabel{AsymptoticFormulaOfMCFOp}
\proof It suffices to prove this relation at $t=0$. However, by hypothesis,
$$
\Psi(\dot{M}(0)) = \dot{\gamma(0)} = -\nabla\opI[u](T_0) = -\Psi(\nabla(\opI[u]\circ\opC)(\opId)).
$$
The result now follows by Lemmas \procref{FirstComponentOfMCFOp} and \procref{SecondComponentOfMCFOp}.\qed
\medskip
\noindent Finally, let $C^{1,\alpha}_{\opin,\opK}(\Bbb{T}\times\Bbb{R})$ and $C^{1,\alpha}_{\opin,\opK^\perp}(\Bbb{T}\times\Bbb{R})$ be the subspaces of $C^{1,\alpha}_\opin(\Bbb{T}\times\Bbb{R})$ consisting of those functions $f:\Bbb{T}\times\Bbb{R}\rightarrow\Bbb{R}$ such that, for all $t$, $f(\cdot,\cdot,t)\in\opK$ and $f(\cdot,\cdot,t)\in\opK^\perp$ respectively.
\proclaim{Theorem \nextprocno}
\noindent If $\opI[u]$ is of Morse-Smale type, then, upon reducing $\epsilon$ if necessary, there exist smooth functions $\eta:]-\epsilon,\epsilon[\rightarrow C^{1,\alpha}(\Bbb{R},\frak{k})$ and $f:]-\epsilon,\epsilon[\rightarrow C^{1,\alpha}_{\opin,\opK^\perp}(\Bbb{T}\times\Bbb{R},]-\pi/4,\pi/4[)$ such that, for all $s$,
$$
\opMCF[\eta(s),f(s),s] = 0.
$$
In particular, for all $s\in]-\epsilon,\epsilon[$, $\ope[\eta(s),f(s),s]$ is a complete mean curvature flow of tori with respect to the metric $g_{1+s}$.
\endproclaim
\proclabel{PerturbationTheorem}
\proof For all $s$, define the linear isomorphism $A_s:C^{0,\alpha}(\Bbb{R},\frak{k})\rightarrow C^{0,\alpha}_{\opin,\opK}(\Bbb{T}\times\Bbb{R})$ by
$$
(A_s\eta)(t) := \Psi(\eta(st)).
$$
Observe that
$$\eqalign{
\|A_s\| &= \opO(\opMax(1,s^\alpha))\ \text{and}\cr
\|A_s^{-1}\| &= \opO(\opMax(1,s^{-\alpha})).\cr}
$$
Define $B:C^{1,\alpha}(\Bbb{R},\frak{k})\rightarrow C^{0,\alpha}(\Bbb{R},\frak{k})$ by
$$
(B\eta)(t) = \frac{\partial\eta}{\partial t}(t) + M(t)^{-1}\opHess(\opI[u]\circ\opC)(M(t))\tau^L\eta(t)(M(t)).
$$
Since $\opI[u]$ is of Morse-Smale type, $B$ is Fredholm and surjective, and thus has a left inverse $\tilde{B}$. Finally, by the classical theory of parabolic operators (see \cite{Taylor}), the operator
$$
C:=\frac{\partial}{\partial t} - (\Delta + 2)
$$
defines a linear isomorphism from $C^{1,\alpha}_{\opin,K^\perp}(\Bbb{T}\times\Bbb{R})$ into $C^{0,\alpha}_{\opin,K^\perp}(\Bbb{T}\times\Bbb{R})$. It follows from the above that the operator
$$
D_s(\eta,f) := A_s B\eta + Cf,
$$
defines a surjective, Fredholm map from $C^{1,\alpha}(\Bbb{R},\frak{k})\oplus C^{1,\alpha}_{\opin,K^\perp}(\Bbb{T}\times\Bbb{R})$ into $C^{0,\alpha}_{\opin}(\Bbb{T}\times\Bbb{R})$ whose left inverse $\tilde{D}_s$ satisfies
$$
\|\tilde{D}_s\| = \opO(\opMax(1,s^{-\alpha})).
$$
However, by Lemma \procref{AsymptoticFormulaOfMCFOp},
$$
\frac{1}{s}\opMCF[\eta,f,s] = D_s(\eta,f) + \opO(s + \|\eta\|_{C^{1,\alpha}}^2),
$$
and the result now follows by the inverse function theorem.\qed
\bigskip
\inappendicestrue
\global\headno=0
\goodbreak
\newhead{Function spaces}[HoelderSpaces]
\newsubhead{H\"older spaces}[HoelderSpaces] We recall the definitions of H\"older space theory that are used in this paper. Let $E$ be a finite-dimensional, normed vector space. Let $M$ be a compact manifold which, for convenience, we assume to be locally isometric to $\Bbb{R}^m$. For all $k\in\Bbb{N}$, let $C^k(M,E)$ be the space of $k$-times continuously differentiable $E$-valued functions over $M$. For $\alpha\in[0,1]$, we define the $\alpha$-H\"older seminorm over $C^0(M,E)$ by
$$
[f]_\alpha := \msup_{x\neq y\in M}\frac{\|f(x)-f(y)\|}{d(x,y)^\alpha}.\eqnum{\nexteqnno[AlphaHoelderSeminorm]}
$$
For all $(k,\alpha)$, we define the $(k,\alpha)$-H\"older norm over $C^k(M,E)$ by
$$
\|f\|_{C^{k,\alpha}} := \sum_{i=0}^k\|D^if\|_{C^0} + [D^kf]_\alpha.\eqnum{\nexteqnno[HoelderNorm]}
$$
We denote by $C^{k,\alpha}(M,E)$ the space of all functions $f\in C^k(M,E)$ such that
$$
\|f\|_{C^{k,\alpha}} < \infty.
$$
We call $C^{k,\alpha}(M,E)$ the H\"older space of order $(k,\alpha)$ over $M$.
\newsubhead{Anisotropic H\"older spaces}[AnisotropicHoelderSpaces]
We consider now function spaces defined over the cartesian product $M\times\Bbb{R}$ where the $M$ and $\Bbb{R}$ components are respectively identified with space and time. Such function spaces are used to study parabolic operators. For $\alpha\in[0,1]$, we define the spatial and temporal H\"older seminorms over $C^0(M\times\Bbb{R},E)$ by
$$\eqalign{
[f]_{x,\alpha} &:= \msup_{x\neq y\in M,t\in\Bbb{R}}\frac{\|f(x,t)-f(y,t)\|}{d(x,y)^\alpha}\ \text{and}\cr
[f]_{t,\alpha} &:= \msup_{x\in M,s\neq t\in\Bbb{R}}\frac{\|f(x,s)-f(x,t)\|}{\left|s-t\right|^\alpha}.\cr}\eqnum{\nexteqnno[SpatialAndTemporalSeminorms]}
$$
For all $k$, we define the anisotropic space $C^k_\opin(M\times\Bbb{R},E)$ to be the space of all functions $f:M\times\Bbb{R}\rightarrow E$ such that $D_x^iD_t^jf$ exists and is continuous whenever $i+2j\leq 2k$, where here $D_x$ and $D_t$ denote respectively the partial derivatives in the spatial and temporal directions. For all $(k,\alpha)$ with $\alpha\in[0,1/2]$, we define the anisotropic $(k,\alpha)$-H\"older norm over $C^k_\opin(M\times\Bbb{R},E)$ by
$$
C^{k,\alpha}_\opin(M\times\Bbb{R},E) := \sum_{i+2j\leq 2k}\|D_x^iD_t^jf\| + \sum_{i+2j=2k}[D_x^iD_t^jf]_{t,\alpha} + \sum_{i+2j=2k}[D_x^iD_t^jf]_{x,2\alpha}.\eqnum{\nexteqnno[AnisotropicHoelderNorm]}
$$
We denote by $C^{k,\alpha}_\opin(M\times\Bbb{R},E)$ the space of all functions $f\in C^k_\opin(M\times\Bbb{R},E)$ such that
$$
\|f\|_{C^{k,\alpha}_\opin} < \infty.
$$
We call $C^{k,\alpha}_\opin(M\times\Bbb{R},E)$ the anisotropic H\"older space of order $(k,\alpha)$ over $M\times\Bbb{R}$.
\newhead{Bibliography}[Bibliography]
{\leftskip = 5ex \parindent = -5ex
\leavevmode\hbox to 4ex{\hfil \cite{Freire}}\hskip 1ex{Alikakos N. D., Freire A., The normalized mean curvature flow for a small bubble in the riemannian manifold, {\sl J. Diff. Geom.}, {\bf 64}, (2003) 247--303}
\medskip
\leavevmode\hbox to 4ex{\hfil \cite{Andrews}}\hskip 1ex{Andrews B., Non-collapsing in mean-convex mean curvature flow, {\sl Geometry \& Topology}, {\bf 16}, (2012), 1413--1418}
\medskip
\leavevmode\hbox to 4ex{\hfil \cite{Brakke}}\hskip 1ex{Brakke K. A., {\sl The motion of a surface by its mean curvature}, Princeton University Press, (1978)}
\medskip
\leavevmode\hbox to 4ex{\hfil \cite{ColdMinII}}\hskip 1ex{Colding T. H., Minicozzi W. P., {\sl A course in minimal surfaces}, Graduate Studies in Mathematics, {\bf 121}, American Mathematical Society, (2011)}
\medskip
\leavevmode\hbox to 4ex{\hfil \cite{ColdMin}}\hskip 1ex{Colding T. H., Ilmanen T., Minicozzi W. P., Rigidity of generic singularities of mean curvature flow, {\sl Publ. Math. Inst. Hautes \'Etudes Sci.}, {\bf 121}, (2015), 363--382}
\medskip
\leavevmode\hbox to 4ex{\hfil \cite{GuilleminPollack}}\hskip 1ex{Guillemin V., Pollack A., {\sl Differential Topology}, American Mathematical Society (Reprint), (2010)}
\medskip
\leavevmode\hbox to 4ex{\hfil \cite{Hamilton}}\hskip 1ex{Hamilton R., Harnack estimate for the mean curvature flow, {\sl J. Diff. Geom}, {\bf 41}, no. 1, (1995), 215--226}
\medskip
\leavevmode\hbox to 4ex{\hfil \cite{Helgason}}\hskip 1ex{Helgason S., {\sl Differential geometry, Lie groups and symmetric spaces}, Academic Press, (1978)}
\medskip
\leavevmode\hbox to 4ex{\hfil \cite{Huisken}}\hskip 1ex{Huisken G., Flow by mean curvature of convex surfaces into spheres, {\sl J. Diff. Geom.}, {\bf 20}, no. 1, (1984), 237--266}
\medskip
\leavevmode\hbox to 4ex{\hfil \cite{Milnor}}\hskip 1ex{Milnor J. W., {\sl Topology from the differential viewpoint}, The University Press of Virginia and Charlottesville, (1965)}
\medskip
\leavevmode\hbox to 4ex{\hfil \cite{MramorPayne}}\hskip 1ex{Mramor A., Payne A., Ancient and Eternal Solutions to Mean Curvature Flow from Minimal Surfaces, arXiv:1904.0843}
\medskip
\leavevmode\hbox to 4ex{\hfil \cite{Salamon}}\hskip 1ex{Robbin J., Salamon D., The spectral flow and the Maslov index, {\sl Bull. London Math. Soc.}, {\bf 27}, no. 1, (1995), 1--33}
\medskip
\leavevmode\hbox to 4ex{\hfil \cite{SmiRos}}\hskip 1ex{Rosenberg H., Smith G, Degree theory of immersed hypersurfaces, to appear in {\sl Mem. Amer. Math. Soc.}, arXiv:1010.1879}
\medskip
\leavevmode\hbox to 4ex{\hfil \cite{Schwarz}}\hskip 1ex{Schwarz M., {\sl Morse Homology}, Birkh\"auser Verlag, (1993)}
\medskip
\leavevmode\hbox to 4ex{\hfil \cite{Smale}}\hskip 1ex{Smale S., An Infinite Dimensional Version of Sard's Theorem, {\sl Amer. J. Math.}, {\bf 87}, no. 4, (1965), 861--866}
\medskip
\leavevmode\hbox to 4ex{\hfil \cite{Taylor}}\hskip 1ex{Taylor M. E., {\sl Partial Differential Equations I}, Applied Mathematics Sciences, {\bf 115}, Springer Verlag, (2011)}
\medskip
\leavevmode\hbox to 4ex{\hfil \cite{White}}\hskip 1ex{White B., The space of minimal submanifolds for varying riemannian metrics, {\sl Ind. Univ. Math. J.}, {\bf 40}, no. 1, (1991), 161--200}
\medskip
\leavevmode\hbox to 4ex{\hfil \cite{WhiteII}}\hskip 1ex{White B., The nature of singularities in mean curvature flow of mean-convex sets, {\sl J. Amer. Math. Soc}, {\bf 16}, no. 1, (2003), 123--138}
\par}
%
%
%
%
\enddocument

%% file: preamble.tex
%
%
%
%
\let\myfrac=\frac%
\input eplain %
\let\frac=\myfrac%
\let\myfootnote=\footnote%
\input amstex \input epsf %
\let\footnote=\myfootnote%
%
%
\loadeufm\loadmsam\loadmsbm\message{symbol names}\UseAMSsymbols\message{,}%
%
\font\myfontdefault=cmr10%
\newif\ifmakebiblio%
\newif\ifinappendices%
\newif\ifundefinedreferences%
\newif\ifchangedreferences%
\makebibliofalse%
\undefinedreferencesfalse%
\changedreferencesfalse%
%
%
%
%
%
\def\setcatcodes{\catcode`\!=0 \catcode`\\=11}%
{\global\let\noe=\noexpand%
\catcode`\@=11 \catcode`\_=11 \setcatcodes%
!global!def!_@@internal@@makeref#1{%
!global!expandafter!def!csname #1ref!endcsname##1{%
!csname _@#1@##1!endcsname%
!expandafter!ifx!csname _@#1@##1!endcsname!relax%
    !write16{#1 ##1 not defined - run saving references}%
    !undefinedreferencestrue%
!fi}}%
!global!def!_@@internal@@makelabel#1{%
!global!expandafter!def!csname #1label!endcsname##1{%
!edef!temptoken{!csname #1info!endcsname}%
!ifloadreferences%
!expandafter!ifx!csname _@#1@##1!endcsname!relax%
!write16{#1 ##1 not hitherto defined - rerun saving references}%
!changedreferencestrue%
!else%
!expandafter!ifx!csname _@#1@##1!endcsname!temptoken%
!else%
!write16{#1 ##1 reference has changed - rerun saving references}%
!changedreferencestrue%
!fi%
!fi%
!else%
!expandafter!edef!csname _@#1@##1!endcsname{!temptoken}%
!edef!textoutput{!write!references{\global\def\_@#1@##1{!temptoken}}}%
!textoutput%
!fi}}%
!global!def!makecounter#1{!_@@internal@@makelabel{#1}!_@@internal@@makeref{#1}}%
!unsetcatcodes%
}
%
%
%
%
%
\def\turnintolatin#1{\ifcase #1 _\or i\or ii\or iii\or iv\or v\or vi\or vii\or viii\or ix\or x\or xi\or xii\or xiii\or xiv\or xv\or xvi\or xvii\or xviii\or xix\or xx\or xxi\or xxii\or xxiii\or xxiv\or xxv\or xxvi\fi}%
\def\alphanum#1{\ifcase #1 _\or A\or B\or C\or D\or E\or F\or G\or H\or I\or J\or K\or L\or M\or N\or O\or P\or Q\or R\or S\or T\or U\or V\or W\or X\or Y\or Z\fi}%
\newwrite\references%
\ifloadreferences{\catcode`\@=11 \catcode`\_=11 \input references.tex }%
\else{\openout\references=references.tex }%
\fi%
%
%
\newcount\headno%
\global\headno=0%
\def\headinfo{\ifinappendices\alphanum\headno\else\the\headno\fi}%
\def\nextheadno{\global\advance\headno by 1 \global\subheadno=0 \global\procno=0 \global\eqnno=0 \headinfo}%
\makecounter{head}%
%
%
\newcount\subheadno%
\global\subheadno=0%
\def\subheadinfo{\headinfo.\the\subheadno}%
\def\nextsubheadno{\global\advance\subheadno by 1 \global\procno=0 \subheadinfo}%
\makecounter{subhead}%
%
%
\newcount\procno%
\global\procno=0%
\def\procinfo{\subheadinfo.\the\procno}%
\def\nextprocno{\global\advance\procno by 1 \procinfo}%
\makecounter{proc}%
%
%
\newcount\figno%
\global\figno=0%
\def\figinfo{\subheadinfo.\the\figno}%
\def\nextfigno{\global\advance\figno by 1 \figinfo}%
\makecounter{fig}%
%
%
\newcount\eqnno%
\global\eqnno=0%
\def\eqninfo{\text{{\rm (\headinfo.\the\eqnno)}}}%
\def\nexteqnno[#1]{\global\advance\eqnno by 1 \eqninfo\hbox{\eqnlabel{#1}}}%
\makecounter{eqn}%
%
%
%
%
%
\def\gobbleeight#1#2#3#4#5#6#7#8{}%
\newcount\citationno%
\global\citationno=0%
\def\citationinfo{\the\citationno}%
\makecounter{citation}%
\newwrite\biblio%
\def\newref#1#2{%
\def\temptext{#2}%
\edef\bibliotextoutput{\expandafter\gobbleeight\meaning\temptext}%
\global\advance\citationno by 1\citationlabel{#1}%
\ifmakebiblio%
    \edef\fileoutput{\write\biblio{\noindent\hbox to 0pt{\hss$[\the\citationno]$}\hskip 0.2em\bibliotextoutput\medskip}}%
    \fileoutput%
\fi}%
\def\cite#1{%
$[\citationref{#1}]$%
\ifmakebiblio%
    \edef\fileoutput{\write\biblio{#1}}%
    \fileoutput%
\fi%
}%
%
%
%
%
\let\mypar=\par%
\edef\Pagetitle={Blank}\headline={\hfil\Pagetitle\hfil}%
\edef\Pagefooter={Blank}\footline={\hfil\Pagefooter\hfil}%
%
%
\newcount\showpagenumflag%
\global\showpagenumflag=0 %
\def\nextoddpage%
{\newpage\ifodd\pageno%
\else\global\showpagenumflag=0 %
\null\vfil\eject%
\global\showpagenumflag=1 %
\fi}%
%
%
\font\headfont=cmb12%
\def\newhead#1[#2]%
{\ifhmode\mypar\fi%
\ifnum\headno=0 \else\goodbreak\bigskip\fi%
{\headfont\noindent\nextheadno\ - #1.}\headlabel{#2}%
\nobreak\medskip}%
%
%
\def\newsubhead#1[#2]%
{\ifhmode\mypar\fi%
\ifnum\subheadno=0 \else\goodbreak\medskip\fi%
{\bf\noindent\nextsubheadno\ - #1.\ }\subheadlabel{#2}}%
%
%
\newif\ifinproclaim%
\global\inproclaimfalse%
\def\proclaim#1{%
\goodbreak\medskip
\bgroup\inproclaimtrue%
\noindent{\bf #1}%
\nobreak\medskip\sl}%
\def\noskipproclaim#1{%
\goodbreak\medskip%
\bgroup\inproclaimtrue%
\noindent{\bf #1}\nobreak\sl}%
\def\endproclaim{\mypar\egroup\nobreak\medskip\ignorespaces}%
%
%
%
\newcount\xpos\newcount\ypos
\def\makelabelgrid{%
\xpos=-5 \ypos=-5 %
\loop\ifnum\xpos<6 %
{\loop\ifnum\ypos<6 %
\def\labeltext{x}%
\ifnum\xpos=0\def\labeltext{+}\fi%
\ifnum\ypos=0\def\labeltext{+}\fi%
\placelabel[\xpos][\ypos]{\labeltext}%
\advance\ypos by 1 %
\repeat}%
\advance\xpos by 1 %
\repeat}%
\def\placelabel[#1][#2]#3{{%
\setbox10=\hbox{\raise #2cm \hbox{\hskip #1cm #3}}%
\ht10=0pt \dp10=0pt \wd10=0pt \box10}}%
%
%
%
%
\def\myitem#1{\noindent\hbox to .5cm{\hfill#1\hss}}%
%
%
%
%
%
%
%
%
%
\font\sansseriften=cmss10%
\font\sansserifseven=cmss7%
\font\sansseriffive=cmss5%
\newfam\sansseriffam%
\textfont\sansseriffam=\sansseriften%
\scriptfont\sansseriffam=\sansserifseven%
\scriptscriptfont\sansseriffam=\sansseriffive%
\def\mathsf{\fam\sansseriffam}%
%
%
%
\font\boldten=cmb10%
\font\boldseven=cmb7%
\font\boldfive=cmb5%
\newfam\mathboldfam%
\textfont\mathboldfam=\boldten%
\scriptfont\mathboldfam=\boldseven%
\scriptscriptfont\mathboldfam=\boldfive%
\def\mathbf{\fam\mathboldfam}%
%
%
%
\font\mycmmiten=cmmi10%
\font\mycmmiseven=cmmi7%
\font\mycmmifive=cmmi5%
\newfam\mycmmifam%
\textfont\mycmmifam=\mycmmiten%
\scriptfont\mycmmifam=\mycmmiseven%
\scriptscriptfont\mycmmifam=\mycmmifive%
\def\hexa#1{\ifcase #1 0\or 1\or 2\or 3\or 4\or 5\or 6\or 7\or 8\or 9\or A\or B\or C\or D\or E\or F\fi}%
\mathchardef\mathi="7\hexa\mycmmifam7B%
\mathchardef\mathj="7\hexa\mycmmifam7C%
%
%
\font\mymsbmten=msbm10 at 8pt%
\font\mymsbmseven=msbm7 at 5.6pt
\font\mymsbmfive=msbm5 at 4pt%
\newfam\mymsbmfam%
\textfont\mymsbmfam=\mymsbmten%
\scriptfont\mymsbmfam=\mymsbmseven%
\scriptscriptfont\mymsbmfam=\mymsbmfive%
\mathchardef\mybeth="7\hexa\mymsbmfam69%
\mathchardef\mygimmel="7\hexa\mymsbmfam6A%
\mathchardef\mydaleth="7\hexa\mymsbmfam6B%
%
%
%
%
\def\proof{{\noindent\bf Proof:\ }}%
\def\remark{{\noindent\bf Remark:\ }}%
\def\qed{~$\square$}%
\def\makeop#1{\global\expandafter\def\csname op#1\endcsname{{\text{#1}}}}%
\def\makeopsmall#1{\global\expandafter\def\csname op#1\endcsname{{\text{\lowercase{#1}}}}}%
%
%
\def\munion{\mathop{\cup}}%
\def\minter{\mathop{\cap}}%
%
%
\makeop{Ext}%
\makeop{Int}%
\makeop{Dist}%
\makeop{Diam}%
\makeop{Length}%
%
%
%
%
%
\def\mlim{\mathop{{\text{Lim}}}}%
\def\mliminf{\mathop{{\text{LimInf}}}}%
\def\msup{\mathop{{\text{Sup}}}}%
%
%
%
\makeop{Dim}%
\makeop{Ker}%
\makeop{Coker}%
\makeop{Tr}%
\makeop{Adj}%
\makeop{Det}%
\makeop{End}%
\makeop{Lin}%
\makeop{Symm}%
\makeop{Mult}%
%
%
\makeop{dx}%
\makeop{dy}%
\makeop{dz}%
\makeop{dt}%
\makeop{dVol}%
\makeop{dArea}%
\makeop{Supp}%
\makeop{Hess}%
\makeop{Lip}%
%
%
\makeop{Re}%
\makeop{Im}%
\makeop{Arg}%
\makeop{Log}%
\makeop{Exp}%
%
%
\makeopsmall{Cos}%
\makeopsmall{Sin}%
\makeopsmall{Tan}%
\makeopsmall{Sec}%
\makeopsmall{Cosec}%
\makeopsmall{Cot}%
\makeopsmall{ArcCos}%
\makeopsmall{ArcSin}%
\makeopsmall{ArcTan}%
\makeopsmall{ArcSec}%
\makeopsmall{ArcCosec}%
\makeopsmall{ArcCot}%
%
%
\makeopsmall{Cosh}%
\makeopsmall{Sinh}%
\makeopsmall{Tanh}%
\makeopsmall{ArcCosh}%
\makeopsmall{ArcSinh}%
\makeopsmall{ArcTanh}%
%
%
\makeop{Vol}%
\makeop{Area}%
\makeop{Riem}%
\makeop{Ric}%
\makeop{Scal}%
\makeop{Euc}%
\makeop{Imm}%
\makeop{Emb}%
%
%
\makeop{Id}%
\makeop{Ad}%
\makeop{O}%
\makeop{SO}%
\makeop{SL}%
\makeop{GL}%
\makeop{Conf}%
\makeop{Homeo}%
\makeop{Diff}%
\makeop{Isom}%
%
%
\makeop{Ind}%
\makeop{Sig}%
\makeop{Spec}%
%
%
\makeop{Conv}%
\makeop{Max}%
\makeop{Min}%
\makeop{Mod}%
\makeop{Deg}%
\makeop{loc}%
%
%
%
%
%
%
%
%
%
%
%
%
%

%% file: references.tex
\global\def\_@citation@Freire{1}
\global\def\_@citation@Andrews{2}
\global\def\_@citation@Brakke{3}
\global\def\_@citation@ColdMinII{4}
\global\def\_@citation@ColdMin{5}
\global\def\_@citation@GuilleminPollack{6}
\global\def\_@citation@Hamilton{7}
\global\def\_@citation@Helgason{8}
\global\def\_@citation@Huisken{9}
\global\def\_@citation@Milnor{10}
\global\def\_@citation@MramorPayne{11}
\global\def\_@citation@Salamon{12}
\global\def\_@citation@SmiRos{13}
\global\def\_@citation@Schwarz{14}
\global\def\_@citation@Smale{15}
\global\def\_@citation@Taylor{16}
\global\def\_@citation@White{17}
\global\def\_@citation@WhiteII{18}
\global\def\_@head@Introduction{1}
\global\def\_@subhead@EternalMeanCurvatureFlows{1.1}
\global\def\_@eqn@MeanCurvatureFlowEquation{\relax \unhbox \voidb@x \hbox {{\relax \tenrm (1.1)}}}
\global\def\_@proc@MainTheoremI{1.1.3}
\global\def\_@subhead@TheMorseSmaleProperty{1.2}
\global\def\_@eqn@StandardCliffordTorusIntro{\relax \unhbox \voidb@x \hbox {{\relax \tenrm (1.2)}}}
\global\def\_@eqn@DefinitionOfFunctionalI{\relax \unhbox \voidb@x \hbox {{\relax \tenrm (1.3)}}}
\global\def\_@proc@MainTheoremII{1.2.1}
\global\def\_@subhead@PerturbingTheMorseComplex{1.3}
\global\def\_@eqn@DefinitionOfMorseChainGroup{\relax \unhbox \voidb@x \hbox {{\relax \tenrm (1.4)}}}
\global\def\_@proc@MainTheoremIII{1.3.1}
\global\def\_@head@TheSpaceOfCliffordTori{2}
\global\def\_@subhead@TheDiffGeomOfEAndCL{2.1}
\global\def\_@eqn@StandardCliffordTorus{\relax \unhbox \voidb@x \hbox {{\relax \tenrm (2.1)}}}
\global\def\_@eqn@StandardParametrisation{\relax \unhbox \voidb@x \hbox {{\relax \tenrm (2.2)}}}
\global\def\_@eqn@FermiParametrisation{\relax \unhbox \voidb@x \hbox {{\relax \tenrm (2.3)}}}
\global\def\_@eqn@FermiMetric{\relax \unhbox \voidb@x \hbox {{\relax \tenrm (2.4)}}}
\global\def\_@eqn@DiffeoOfOWithCL{\relax \unhbox \voidb@x \hbox {{\relax \tenrm (2.5)}}}
\global\def\_@eqn@StabiliserOfTZero{\relax \unhbox \voidb@x \hbox {{\relax \tenrm (2.6)}}}
\global\def\_@eqn@DefinitionOfHatE{\relax \unhbox \voidb@x \hbox {{\relax \tenrm (2.7)}}}
\global\def\_@eqn@DefinitionOfE{\relax \unhbox \voidb@x \hbox {{\relax \tenrm (2.8)}}}
\global\def\_@eqn@JacobiOperator{\relax \unhbox \voidb@x \hbox {{\relax \tenrm (2.9)}}}
\global\def\_@eqn@TangentSpaceOfCL{\relax \unhbox \voidb@x \hbox {{\relax \tenrm (2.10)}}}
\global\def\_@subhead@TheSymmetricSpaceStructuresOfSAndCL{2.2}
\global\def\_@eqn@DefinitionOfB{\relax \unhbox \voidb@x \hbox {{\relax \tenrm (2.11)}}}
\global\def\_@eqn@DefinitionOfHAndK{\relax \unhbox \voidb@x \hbox {{\relax \tenrm (2.12)}}}
\global\def\_@eqn@DefinitionOfXiA{\relax \unhbox \voidb@x \hbox {{\relax \tenrm (2.13)}}}
\global\def\_@eqn@MetricOverEnd{\relax \unhbox \voidb@x \hbox {{\relax \tenrm (2.14)}}}
\global\def\_@eqn@FormulaForPhiA{\relax \unhbox \voidb@x \hbox {{\relax \tenrm (2.15)}}}
\global\def\_@subhead@KillingFieldsOverCL{2.3}
\global\def\_@eqn@FlowsAndExponentialFlows{\relax \unhbox \voidb@x \hbox {{\relax \tenrm (2.16)}}}
\global\def\_@eqn@DefinitionOfXAndF{\relax \unhbox \voidb@x \hbox {{\relax \tenrm (2.17)}}}
\global\def\_@eqn@ExplicitFormulaForX{\relax \unhbox \voidb@x \hbox {{\relax \tenrm (2.18)}}}
\global\def\_@proc@ExplicitFormulaForX{2.3.1}
\global\def\_@subhead@TheGeometryOfCurvesInCL{2.4}
\global\def\_@eqn@DefinitionOfGammaInCL{\relax \unhbox \voidb@x \hbox {{\relax \tenrm (2.19)}}}
\global\def\_@eqn@DefinitionOfFtInCL{\relax \unhbox \voidb@x \hbox {{\relax \tenrm (2.20)}}}
\global\def\_@eqn@AsymptoticFormulaForF{\relax \unhbox \voidb@x \hbox {{\relax \tenrm (2.21)}}}
\global\def\_@proc@AsymptoticFormulaForF{2.4.1}
\global\def\_@eqn@ExponentialMapOfCL{\relax \unhbox \voidb@x \hbox {{\relax \tenrm (2.22)}}}
\global\def\_@eqn@FormulaForGamma{\relax \unhbox \voidb@x \hbox {{\relax \tenrm (2.23)}}}
\global\def\_@eqn@DefinitionOfAlphaAndH{\relax \unhbox \voidb@x \hbox {{\relax \tenrm (2.24)}}}
\global\def\_@eqn@AsymptoticFormulaForFConstructingCurves{\relax \unhbox \voidb@x \hbox {{\relax \tenrm (2.25)}}}
\global\def\_@proc@AsymptoticFormulaForFConstructingCurves{2.4.2}
\global\def\_@eqn@DerivativesOfHAlphaBeta{\relax \unhbox \voidb@x \hbox {{\relax \tenrm (2.26)}}}
\global\def\_@proc@DerivativesOfHAlphaBeta{2.4.3}
\global\def\_@eqn@DerivativesOfHAlpha{\relax \unhbox \voidb@x \hbox {{\relax \tenrm (2.27)}}}
\global\def\_@subhead@ZeroSetsOfTheKernel{2.5}
\global\def\_@eqn@DefinitionOfZA{\relax \unhbox \voidb@x \hbox {{\relax \tenrm (2.28)}}}
\global\def\_@eqn@DefinitionOfJ{\relax \unhbox \voidb@x \hbox {{\relax \tenrm (2.29)}}}
\global\def\_@eqn@RotationsYieldTranslations{\relax \unhbox \voidb@x \hbox {{\relax \tenrm (2.30)}}}
\global\def\_@eqn@DefinitionOfPortionsOfSigma{\relax \unhbox \voidb@x \hbox {{\relax \tenrm (2.31)}}}
\global\def\_@eqn@DefinitionOfTheta{\relax \unhbox \voidb@x \hbox {{\relax \tenrm (2.32)}}}
\global\def\_@eqn@DefinitionOfThetaPM{\relax \unhbox \voidb@x \hbox {{\relax \tenrm (2.33)}}}
\global\def\_@eqn@DefinitionOfD{\relax \unhbox \voidb@x \hbox {{\relax \tenrm (2.34)}}}
\global\def\_@proc@MatrixDecompositionLemma{2.5.1}
\global\def\_@eqn@EquivarianceOfTheta{\relax \unhbox \voidb@x \hbox {{\relax \tenrm (2.35)}}}
\global\def\_@eqn@StandardSingularMatrix{\relax \unhbox \voidb@x \hbox {{\relax \tenrm (2.36)}}}
\global\def\_@eqn@SingularZeroSet{\relax \unhbox \voidb@x \hbox {{\relax \tenrm (2.37)}}}
\global\def\_@eqn@StandardGenericMatrixII{\relax \unhbox \voidb@x \hbox {{\relax \tenrm (2.38)}}}
\global\def\_@eqn@ZeroSetSubstitution{\relax \unhbox \voidb@x \hbox {{\relax \tenrm (2.39)}}}
\global\def\_@eqn@NonSingularNonSpecialZeroSetII{\relax \unhbox \voidb@x \hbox {{\relax \tenrm (2.40)}}}
\global\def\_@eqn@SecondDerivative{\relax \unhbox \voidb@x \hbox {{\relax \tenrm (2.41)}}}
\global\def\_@proc@TheSecondDerivativeIsNonVanishing{2.5.2}
\global\def\_@head@TheMorseSmaleProperty{3}
\global\def\_@subhead@ConstructingFunctionsOfMorseType{3.1}
\global\def\_@eqn@DefinitionOfI{\relax \unhbox \voidb@x \hbox {{\relax \tenrm (3.1)}}}
\global\def\_@proc@TheMorseProperty{3.1.1}
\global\def\_@eqn@FirstJetOfU{\relax \unhbox \voidb@x \hbox {{\relax \tenrm (3.2)}}}
\global\def\_@proc@SardSmaleTheorem{3.1.2}
\global\def\_@proc@PreSurjectivity{3.1.3}
\global\def\_@proc@FirstSubmersionLemma{3.1.4}
\global\def\_@subhead@TheSmaleProperty{3.2}
\global\def\_@eqn@GradientFlowOperator{\relax \unhbox \voidb@x \hbox {{\relax \tenrm (3.3)}}}
\global\def\_@eqn@DefinitionOfM{\relax \unhbox \voidb@x \hbox {{\relax \tenrm (3.4)}}}
\global\def\_@eqn@DOneGradientFlowOperator{\relax \unhbox \voidb@x \hbox {{\relax \tenrm (3.5)}}}
\global\def\_@proc@TheMorseSmaleProperty{3.2.1}
\global\def\_@eqn@DTwoGradentFlowOperator{\relax \unhbox \voidb@x \hbox {{\relax \tenrm (3.6)}}}
\global\def\_@proc@SurjectivityMorseSmale{3.2.2}
\global\def\_@eqn@DualDOneGradientFlowOperator{\relax \unhbox \voidb@x \hbox {{\relax \tenrm (3.7)}}}
\global\def\_@eqn@PairingToProveSurjectivity{\relax \unhbox \voidb@x \hbox {{\relax \tenrm (3.8)}}}
\global\def\_@proc@SurjectivityMorseSmaleII{3.2.3}
\global\def\_@proc@GIsOrthogonalToFlow{3.2.4}
\global\def\_@subhead@TowardsTheMorseSmaleProperty{3.3}
\global\def\_@eqn@HigherOrderJets{\relax \unhbox \voidb@x \hbox {{\relax \tenrm (3.9)}}}
\global\def\_@proc@SecondSubmersionLemma{3.3.1}
\global\def\_@proc@FirstTransversalProperty{3.3.2}
\global\def\_@proc@SecondTransversalProperty{3.3.3}
\global\def\_@proc@ThirdTransversalProperty{3.3.4}
\global\def\_@subhead@ConstructingThePerturbation{3.4}
\global\def\_@eqn@DefinitionOfUThree{\relax \unhbox \voidb@x \hbox {{\relax \tenrm (3.10)}}}
\global\def\_@eqn@PositivityOfG{\relax \unhbox \voidb@x \hbox {{\relax \tenrm (3.11)}}}
\global\def\_@eqn@NonVanishingOfSecondDerivative{\relax \unhbox \voidb@x \hbox {{\relax \tenrm (3.12)}}}
\global\def\_@eqn@FirstComponentOfV{\relax \unhbox \voidb@x \hbox {{\relax \tenrm (3.13)}}}
\global\def\_@eqn@SecondComponentOfV{\relax \unhbox \voidb@x \hbox {{\relax \tenrm (3.14)}}}
\global\def\_@eqn@DefinitionOfTildeVDelta{\relax \unhbox \voidb@x \hbox {{\relax \tenrm (3.15)}}}
\global\def\_@eqn@DefinitionOfVDelta{\relax \unhbox \voidb@x \hbox {{\relax \tenrm (3.16)}}}
\global\def\_@eqn@BoundsAwayFromZero{\relax \unhbox \voidb@x \hbox {{\relax \tenrm (3.17)}}}
\global\def\_@proc@BoundsAwayFromZero{3.4.1}
\global\def\_@eqn@ConstructingVPartI{\relax \unhbox \voidb@x \hbox {{\relax \tenrm (3.18)}}}
\global\def\_@proc@ConstructingVPartI{3.4.2}
\global\def\_@eqn@ConstructingVPartII{\relax \unhbox \voidb@x \hbox {{\relax \tenrm (3.19)}}}
\global\def\_@proc@ConstructingVPartII{3.4.3}
\global\def\_@eqn@ConstructingVPartIII{\relax \unhbox \voidb@x \hbox {{\relax \tenrm (3.20)}}}
\global\def\_@proc@ConstructingVPartIII{3.4.4}
\global\def\_@subhead@TechnicalPropertiesOfV{3.5}
\global\def\_@eqn@MorseSmaleAsymptoticsOfG{\relax \unhbox \voidb@x \hbox {{\relax \tenrm (3.21)}}}
\global\def\_@eqn@AsymptoticFormulaForIntegrand{\relax \unhbox \voidb@x \hbox {{\relax \tenrm (3.22)}}}
\global\def\_@eqn@MorseSmaleNormalField{\relax \unhbox \voidb@x \hbox {{\relax \tenrm (3.23)}}}
\global\def\_@eqn@MorseSmaleAreaForm{\relax \unhbox \voidb@x \hbox {{\relax \tenrm (3.24)}}}
\global\def\_@eqn@MorseSmaleAsymptoticsOfIntegrand{\relax \unhbox \voidb@x \hbox {{\relax \tenrm (3.25)}}}
\global\def\_@eqn@PropertiesOfTildePhi{\relax \unhbox \voidb@x \hbox {{\relax \tenrm (3.26)}}}
\global\def\_@eqn@MorseSmaleAsymptoticsOfF{\relax \unhbox \voidb@x \hbox {{\relax \tenrm (3.27)}}}
\global\def\_@eqn@MorseSmaleAAndB{\relax \unhbox \voidb@x \hbox {{\relax \tenrm (3.28)}}}
\global\def\_@eqn@EstimateOverEtaInterval{\relax \unhbox \voidb@x \hbox {{\relax \tenrm (3.29)}}}
\global\def\_@proc@EstimateOverEtaInterval{3.5.1}
\global\def\_@eqn@MorseSmaleDerivativeOfF{\relax \unhbox \voidb@x \hbox {{\relax \tenrm (3.30)}}}
\global\def\_@eqn@DefinitionOfTildePhi{\relax \unhbox \voidb@x \hbox {{\relax \tenrm (3.31)}}}
\global\def\_@eqn@FirstRelationOfPhis{\relax \unhbox \voidb@x \hbox {{\relax \tenrm (3.32)}}}
\global\def\_@eqn@MSEstimateI{\relax \unhbox \voidb@x \hbox {{\relax \tenrm (3.33)}}}
\global\def\_@proc@MSEstimateI{3.5.2}
\global\def\_@eqn@MSEstimateII{\relax \unhbox \voidb@x \hbox {{\relax \tenrm (3.34)}}}
\global\def\_@proc@MSEstimateII{3.5.3}
\global\def\_@eqn@MSEstimateIII{\relax \unhbox \voidb@x \hbox {{\relax \tenrm (3.35)}}}
\global\def\_@proc@MSEstimateIII{3.5.4}
\global\def\_@eqn@MSEstimateIVA{\relax \unhbox \voidb@x \hbox {{\relax \tenrm (3.36)}}}
\global\def\_@eqn@MSEstimateIVB{\relax \unhbox \voidb@x \hbox {{\relax \tenrm (3.37)}}}
\global\def\_@proc@MSEstimateIV{3.5.5}
\global\def\_@eqn@MSEstimateV{\relax \unhbox \voidb@x \hbox {{\relax \tenrm (3.38)}}}
\global\def\_@proc@MSEstimateV{3.5.6}
\global\def\_@eqn@MSEstimateVI{\relax \unhbox \voidb@x \hbox {{\relax \tenrm (3.39)}}}
\global\def\_@proc@MSEstimateVI{3.5.7}
\global\def\_@head@PerturbationTheory{4}
\global\def\_@subhead@WhitesConstruction{4.1}
\global\def\_@eqn@DerivativeOfMeanCurvatureFunctional{\relax \unhbox \voidb@x \hbox {{\relax \tenrm (4.1)}}}
\global\def\_@eqn@DefinitionOfFI{\relax \unhbox \voidb@x \hbox {{\relax \tenrm (4.2)}}}
\global\def\_@eqn@DefinitionOfFII{\relax \unhbox \voidb@x \hbox {{\relax \tenrm (4.3)}}}
\global\def\_@eqn@DefinitionOfH{\relax \unhbox \voidb@x \hbox {{\relax \tenrm (4.4)}}}
\global\def\_@eqn@LeftInvariance{\relax \unhbox \voidb@x \hbox {{\relax \tenrm (4.5)}}}
\global\def\_@eqn@RightInvariance{\relax \unhbox \voidb@x \hbox {{\relax \tenrm (4.6)}}}
\global\def\_@subhead@TheDerivativesOfA{4.2}
\global\def\_@eqn@FirstDerivativeOfTildeAWRTG{\relax \unhbox \voidb@x \hbox {{\relax \tenrm (4.7)}}}
\global\def\_@proc@FirstDerivativeOfTildeAII{4.2.1}
\global\def\_@eqn@FirstDerivativeOfTildeA{\relax \unhbox \voidb@x \hbox {{\relax \tenrm (4.8)}}}
\global\def\_@proc@FirstDerivativeOfTildeA{4.2.2}
\global\def\_@eqn@TildeHIsSmall{\relax \unhbox \voidb@x \hbox {{\relax \tenrm (4.9)}}}
\global\def\_@eqn@HessianOfArea{\relax \unhbox \voidb@x \hbox {{\relax \tenrm (4.10)}}}
\global\def\_@proc@HessianOfArea{4.2.3}
\global\def\_@subhead@TheMeanCurvatureFlowOperator{4.3}
\global\def\_@eqn@DerivativeOfFamilyOfMetrics{\relax \unhbox \voidb@x \hbox {{\relax \tenrm (4.11)}}}
\global\def\_@eqn@DefinitioOfLiftOfGamma{\relax \unhbox \voidb@x \hbox {{\relax \tenrm (4.12)}}}
\global\def\_@eqn@DefinitionOfMCF{\relax \unhbox \voidb@x \hbox {{\relax \tenrm (4.13)}}}
\global\def\_@eqn@AsymptoticFormulaForNormalVariation{\relax \unhbox \voidb@x \hbox {{\relax \tenrm (4.14)}}}
\global\def\_@proc@FirstComponentOfMCFOp{4.3.1}
\global\def\_@eqn@AsymptoticFormulaForMeanCurvature{\relax \unhbox \voidb@x \hbox {{\relax \tenrm (4.15)}}}
\global\def\_@proc@SecondComponentOfMCFOp{4.3.2}
\global\def\_@eqn@AsymptoticFormulaForMeanCurvatureFlowOperator{\relax \unhbox \voidb@x \hbox {{\relax \tenrm (4.16)}}}
\global\def\_@proc@AsymptoticFormulaOfMCFOp{4.3.3}
\global\def\_@proc@PerturbationTheorem{4.3.4}
\global\def\_@head@HoelderSpaces{A}
\global\def\_@subhead@HoelderSpaces{A.1}
\global\def\_@eqn@AlphaHoelderSeminorm{\relax \unhbox \voidb@x \hbox {{\relax \tenrm (A.1)}}}
\global\def\_@eqn@HoelderNorm{\relax \unhbox \voidb@x \hbox {{\relax \tenrm (A.2)}}}
\global\def\_@subhead@AnisotropicHoelderSpaces{A.2}
\global\def\_@eqn@SpatialAndTemporalSeminorms{\relax \unhbox \voidb@x \hbox {{\relax \tenrm (A.3)}}}
\global\def\_@eqn@AnisotropicHoelderNorm{\relax \unhbox \voidb@x \hbox {{\relax \tenrm (A.4)}}}
\global\def\_@head@Bibliography{B}